\documentclass{article}
\usepackage{graphicx, pdfsync}
\usepackage[colorlinks=true, citecolor=red]{hyperref}
\usepackage{amsmath,amssymb,amsfonts,MnSymbol, mathrsfs}

\usepackage{algorithm}
\usepackage{xcolor}
\usepackage{listings}
\usepackage[colorlinks=true, citecolor=red]{hyperref}

\usepackage{tikz,pgfplots,bm}
\usetikzlibrary{patterns}

\newtheorem{theorem}{Theorem}
\newtheorem{proposition}[theorem]{Proposition}%
\newtheorem{remark}{Remark}%
\newtheorem{notation}{Notation}

\def \e {{\rm e}}
\def \R{{\mathbb{R}}}
\def \S{{\mathbb{S}}}

\newcommand{\ds}{\displaystyle}
\newcommand{\be}{\begin{equation}}
\newcommand{\ee}{\end{equation}}

\def\p{\partial}
\def\n{\nabla}
\def\nablav{\bm\nabla}
\def\d{\hbox{d}}
\def\One{{\bf 1}}

\def\vx{{\bf x}}

\def\vu{{\bf u}}
\def\omegav{{\bm\omega}}
\def\vn{{\bf n}}

\def\vom{{\bm\omega}}
\def \SS{{\mathbb{S}^2}}

\def\d{{\partial}}
\def\dd{\mathrm{d}}

\def\bR{\mathbf{R}}

\title{{\Huge Radiative Transfer For Variable 3D Atmospheres}
\\
Working Document}
\author{F. Golse, F. Hecht,  O. Pironneau, D. Smets and P.-H. Tournier}

\begin{document}

\maketitle

\begin{abstract}To study the temperature in a gas subjected to electromagnetic radiations, one may use the Radiative Transfer equations coupled with the Navier-Stokes equations.  The problem has 7 dimensions;  however with minimal simplifications it is equivalent  to a small number of integro-differential equations in 3 dimensions. We present the method and a numerical implementation using an ${\mathcal H}$-matrix compression scheme. The result is a very fast: 50K physical points, all directions of radiation and 680 frequencies require less than 5 minutes on an Apple M1 Laptop. The method is capable of handling variable absorptioN and scattering functionS of spatial positions and frequencies.

  The implementation is done using \texttt{htool}\footnote{https://github.com/htool-ddm/htool}, a matrix compression library interfaced with the PDE solver\texttt{freefem++}. Applications to the temperature in the French Chamonix valley is presented at different hours of the day with and without snow / clouds and with a variable absorption taken from the Gemini measurements. The result is precise enough to assert temperature differences due to increased absorption in the vibrational frequency subrange of greenhouse gasses.
\end{abstract}


\section*{Introduction}
Heat transfer with radiative transfer are very important in astronomy, combustion and climatology, to name a few.  For the atmosphere the reader is sent to \cite{clouds},\cite{GOO}, \cite{BOH},  the numerically oriented book \cite{ZDU} and the two mathematically oriented books \cite{CHA} and \cite{FOW}.
Due to the considerable difference of length  scales, modelling at the level of photons is difficult to use for the earth atmosphere in large areas. A much simpler formulation, known as the radiative transfer equations, is based on energy conservation principles of continuum mechanics. 

These were shown to be well posed in \cite{POR} for the Radiative Transfer system coupled with the time dependent heat equation. Existence and uniqueness was proved in \cite{FGOP} and \cite{DIA} for the stationary case.

On the numerical side the problem is hard because the model has 3 space variables, two directional ones, the frequencies and time. 

In 2005, K. Evans and A. Marshak wrote in chapter 4 of \cite{clouds} a  review of the numerical methods available for Radiative Transfer alone. Today, judging from \cite{I3RC}, the situation has not changed: SHDOM (Spherical Harmonic Discrete Ordinate Method) and Monte-Carlo are the two most popular methods. While reviewing the current situation for the radiative transfer equations in \cite{CBOP} we implemented a finite element version of SHDOM and found that the method was incapable, unless a huge number of degree of freedom is used, of giving results with the accuracy needed to differentiate between small variations of the absorption coefficient. 

One may resort to approximations. In nuclear physics and astronomy,  for the numerical simulations in the 1960 a constant  absorption coefficient was used, the so called Grey Model. In climatology, the grey model cannot explain the role of the greenhouse gasses.

   The stratified approximation assumes that the surface which receives the radiations is flat and the source is far.  Thus only one space variable is retained. For such cases, an integral formulation, probably due to Chandrasekhar \cite{CHA}, turns out to be much more precise and also computationally much cheaper. A fixed-point iteration scheme on this nonlinear integral formulation, known in the radiative transfer community as ``iterations on the 
sources" was shown to be monotone in \cite{OP}, a property which seems to have escaped earlier studies.   Finally in \cite{FGOP},and \cite{FGOP3} the method was extended to include the temperature equation of the fluid and also to handle Rayleigh scattering while retaining 
monotonicity. 
In \cite{DIA} F. Golse observed that the integral formulation exists also in 3D. A new existence result was proved by using this formulation and it was numerically tested.  However the computing cost was $O(N^2)$, $N$ being the number of vertices of the triangulation of the physical domain.

The purpose of this article is to present an $O(N\ln N)$ implementation which uses compressed ${\mathcal H}$-matrices and views the computation of integrals as a matrix-vector product in the Finite Element discrete space.  Thence P.H. Tournier, who is a co-author of the library \texttt{htool} for integral equations with ${\mathcal H}$-matrices, wrote the necessary modifications for Radiative Transfer, namely the handling of mixed matrices with one vertex on the boundary and the other one inside the domain.

The programming is a joint work of F. Hecht, O. Pironneau and P.H. Tournier, done using the high level PDE solver \texttt{freefem++} \cite{FF}. 

The method is tested on a $35\times 35$ km area and the first 10km of the atmosphere above it.  The center of the domain is the city of Chamonix in the French Alps, provided by D. Smets who has constructed high precision maps and an automatic triangular mesh generator for many places in the world for the Microsoft flight simulator. Hence this Radiative-Heat-Transfer code can be used easily for any other terrain. Note that such a study would be a challenge for Cartesian meshes because some of the mountains like the Mont-Blanc are above 4000m and some slopes are almost vertical.

The absorption coefficient is a function of space and frequencies. It could also be a function of temperature if another fixed point  loop is added.  No regularity is assumed beside positivity and boundedness.
Similarly, the scattering variable coefficient need only to be in (0,1).
  
Real life data are used; the scattering part is not yet implemented As shown in \cite{FGOP2}, it is not so difficult to do, but it adds complexity to an already complex topic.  
Numerically, the methods is very fast: it can handle a few hundred thousand physical points plus all directions plus 680 frequency points on a PC (and a few millions on a supercomputer) within minutes.

Note that the authors have no competence for climate modelling, so the results are only briefly commented:  the purpose of this study is to validate the numerical tool.

\section{The Mathematical Problem}
To find the temperature $T$ in an incompressible fluid exposed to electromagnetic radiations, it is necessary to solve the Navier-Stokes equations coupled with the Radiative Transfer equations.  It is a system of partial differential equations formulated in terms of the fluid velocity $\vu$, its pressure $p$, its density $\rho$, and temperature $T$; all are  functions of time $t$ and position $\vx$ in the physical domain $\Omega$.  It involves also the light intensity field ${I_\nu}(\vx,t,\omegav,t)$ 
for each frequency $\nu$  in each direction $\omegav$:

\smallskip

Given ${I_\nu},T,\vu,\rho$ at $t=0$, find ${I_\nu},T, \vu,\rho$, s.t. $\{\vx,\omegav,t,\nu\}\in\Omega\times\S_2\times(0,\bar T)\times\R^+$,
\begin{equation}
\begin{aligned}
\label{onea}
&\frac1c\p_t {I_\nu} + \omegav\cdot\nablav {I_\nu}+\rho\tilde\kappa a\left[{I_\nu}-{\frac{1}{4\pi}\int_{\S^2}} p_\nu(\omegav,\omegav'){I_\nu}(\omegav'){d}\omega'\right]
\\&
\hskip6cm = \rho\tilde\kappa(1-a) [B_\nu(T)-{I_\nu}],
\\ 
&\rho (\p_tT+\vu\cdot\n T) -\nabla\cdot(\rho \kappa_T\nabla T )
+ A\nablav\cdot\int_0^\infty\frac1{4\pi}\int_{\S^2}
 {I_\nu}(\omegav)\omegav \dd\omega\dd\nu =0
\\ 
&\p_t\vu+\vu\cdot\n\vu-\frac{\mu_F}\rho\Delta\vu + \frac1\rho\n p={\bf g}(T),\quad \n\cdot\vu=0,
\quad 
\p_t\rho+\n\cdot(\rho\vu)=0,
\end{aligned}
\end{equation}
where $A=4\pi/c_P$, $c_P$ is the specific heat of the medium at constant pressure,  $\SS$ is the unit sphere, $\nabla,\Delta$ are with respect to $\vx$, 
$\ds B_\nu(T)=\frac{2 \hbar \nu^3}{c^2[{\rm e}^\frac{\hbar\nu}{k T}-1]}$, is the Planck function,
  $\hbar,c,k$ are the Planck constant, the speed of light in the medium and  
the Boltzmann constant.
The above holds only if $c_P$ is constant in $\Omega$, which is the case in moderate size domains.

\subsection{Adimensionalization}
Recall that
$
c=2.99\cdot 10^9 m\cdot s^{-1},~  \hbar=6.63\cdot 10^{-34} J.s,~ k=1.38 \cdot 10^{-23} J\cdot K^{-1},
$
\[
\kappa_T=2\cdot 10^{-5} m^2\cdot s^{-1}, \quad c_P=1.1~ J.s\cdot(g.K)^{-1}, \quad \rho \sim \bar\rho:=1.3 g\cdot m^{-3},\hbox{ in air at $20^o$C}. 
\]%
The frequencies of interest are in the order of $\nu_0=10^{14}s^{-1}$, so we define  $B_0=\frac{2\hbar\nu_0^3}{c^2}$, $T_0=\frac{\hbar\nu_0}k$ and work with $I_\nu/B_0$, $T/T_0$ and $\nu/\nu_0$. Then  \eqref{onea} holds with these new unknowns with $\ds B_\nu(T)={\nu^3}/({\rm e}^\frac\nu T -1)$, $c_P=1$ and
\begin{equation}
\label{rescale}
B_0=1.47~J\cdot m^{-2}\quad T_0=4789 K, \quad  \frac A{\bar\rho}= \frac{4\pi B_0\nu_0}{c_P\bar\rho T_0}=2.70\cdot 10^{11}~ m\cdot s^{-1}.
\end{equation}
\begin{notation}
If $\bar\rho$ is a reference density for the gas, $\kappa:=\bar\rho\tilde\kappa$ is   the percentage of light absorbed per unit length. It has the dimension $m^{-1}$, thence, the notation
\[
\kappa(\vx,\nu) := \rho(\vx)\tilde\kappa_\nu(\vx).  
\]
Later we may assume at times that $\tilde\kappa_\nu(\vx)$ is independent of $\vx$ or is the sum of products of $\vx$-functions with $\nu$-functions.
\end{notation}
   The scattering albedo is  $a\in(0,1)$ and $\frac1{4\pi}p_\nu(\omegav,\omegav')$ is the ``probability''  that a ray in direction $\omegav'$ scatters in direction $\omegav$. Both $\kappa$ and $a$ are functions of $\nu$. 
   We further assume that $\nu\mapsto\kappa$ and $\nu\mapsto a$ are continuous positive functions, satisfying, for some positive constants $a_M$ and $\kappa_m<\kappa_M$:
$$
0\le\kappa_m\le\kappa\le\kappa_M\,,\qquad 0\le a\le a_M<1\,,\qquad\forall\nu>0\,.
$$
The viscosity of air $\mu_F$ is $18\cdot 10^{-6}$Pa.s; ${\bf g}(T)$, the Boussinesq term, is a vector valued function of the temperature $T$. 
\smallskip

We require $\Omega$ to be an open bounded subset of $\bR^3$ with $C^1$ boundary.
We denote by $\vn(\vx)$ the outward unit normal to $\Gamma:=\d\Omega$ at $\vx$. 

Boundary conditions are needed, for example, initial values for $I_\nu,T,\vu,\rho$,  $\vu$ given on $\Gamma$ and $\rho$ given at $\vx\in\Gamma$ where $\vu\cdot\vn<0$ and 
\[
I_\nu(\vx,\vom)\!=\!Q_\nu(\vx,\vom)\,,\quad\vom\cdot\vn<0\,,\qquad\frac{\d T}{\d n}\Big\vert _{\Gamma}=0\,.
\]
\subsubsection{Sunlight}
When $Q_\nu(\vx,\vom)$ is due to sunlight, of power $Q^0$ coming from the direction $\omegav_s$, then the surface which receives the light re-emits it with (Lambertian reflexion), 
\[
Q_\nu(\vx,\vom)=Q_\nu(\vx)\cdot[\omegav\cdot\vn(\vx)]_-\hbox{ with }Q_\nu(\vx):=Q^0 [\omegav_s\cdot\vn(\vx)]_+ B_\nu(T_s),
\] 
the Planck function with the Sun's temperature $T_s$.

\subsection{Simplifications}

The angular average radiative intensity plays an important role in this article: 
\[
J_\nu(\vx,t) = \frac1{4\pi}\int_{\S^2} {I_\nu}(\omegav ){d}\omega.
\]
Averaging the first equation in \eqref{onea} and neglecting the time dependent term because $c>>1$, gives :
\begin{equation}\label{oneb}
\nablav\cdot\int_0^\infty{\int_{\S^2}} {I_\nu}(\omegav)\omegav {d}\omega {d}\nu
 =  4\pi\int_0^\infty  \rho\tilde\kappa(1-a)\left( B_\nu(T)-J_\nu\right) {d}\nu
\,,
\end{equation}
which, in turn, leads to a simpler form ( see \eqref{RTHeat}) for the second equation in \eqref{onea}.

\subsubsection{Isotropic Scattering}
Assume isotropic scattering: $p_\nu(\omegav,\omegav')\equiv 1$. Let us neglect the variations of $\rho$ in the diffusion term of the temperature equation. Let us neglect $\frac1c\partial_t I$ because $c$ is very large and let us assume that the Boussinesq term is small.  Then the Navier-Stokes equations are decoupled and can be solved before hand. Then in the Radiative transfer and temperature equations $\rho$, $\vu$  are given. So \eqref{onea} becomes the following system for $I_\nu,J_\nu,T$:
\begin{equation}\label{RTHeat}
\left\{
\begin{aligned}
{}&\vom\cdot\nabla I_\nu+{\rho\tilde\kappa} I_\nu={\rho\tilde\kappa}(1-a)B_\nu(T)+{\rho\tilde\kappa} a J_\nu\,,\quad J_\nu:=\tfrac1{4\pi}\int_{\SS}I_\nu\dd\vom\,,
\\
&\p_tT+\vu\cdot\n T -\nabla\cdot(\kappa_T\nabla T )= A\int_0^\infty{\tilde\kappa}(1-a)(J_\nu-B_\nu(T))\dd\nu\,,
\\
&I_\nu(\vx,\vom)\!=\!Q_\nu(\vx,\vom)\,,\quad\vom\cdot\vn<0\,,\,\,\vx\in\d\Omega\,,\qquad\frac{\d T}{\d n}\Big\vert _{\d\Omega}=0\,,
\quad T\big\vert _{t=0}=T_{0}\,.
\end{aligned}
\right.
\end{equation}
\subsubsection{Small Thermal Diffusion and Convection}
When the convective velocity and $\kappa_T$  are small compared to $A$, the temperature equation reaches a stationary states and simplifies to an integral relation between $J$ and $B_\nu(T)$ and the system becomes:
%
\begin{align}
&\vom\cdot\nabla I_\nu+{\rho\tilde\kappa} I_\nu={\rho\tilde\kappa}(1-a)B_\nu(T)+{\rho\tilde\kappa} a J_\nu\,,\quad J_\nu:=\tfrac1{4\pi}\int_{\SS}I_\nu\dd\vom\,,
\label{RTa}
\\
&\int_0^\infty{\rho\tilde\kappa}(1-a)(J_\nu-B_\nu(T))\dd\nu=0\,,
 \label{RTb}
 \\ &
I_\nu(\vx,\vom) =Q_\nu(\vx,\vom)
\, \hbox{ on }\Sigma:=\{(\vx,\vom)\in\Gamma\times\SS~:~\vom(\vx)\cdot\vn(\vx)<0\,\}\,.
\label{RTc}
\end{align}
The system is 6 dimensional in $\left[\vx=[x,y,z],\vom=[\omega_1,\omega_2,\omega_3~:~|\vom|=1],\nu\right]$. 

\subsubsection{The stratified case}
When $\Omega$ is thin in a direction $Ox$ and the physical data depend mostly on $x$, the equations becomes almost one dimensional in $x$. The domain $\Omega$ is locally $(0,H)\times\R^2$.
Assuming that sunlight at $x=+\infty$ crosses the atmosphere unaffected
 the equations reduce to
\begin{equation}
\begin{aligned}
&\mu\partial_x I_\nu+{\rho\tilde\kappa} I_\nu={\rho\tilde\kappa}(1-a)B_\nu(T)+{\rho\tilde\kappa} a J_\nu\,,\quad J_\nu:=\tfrac1{2}\int_0^1 I_\nu \dd\mu
\label{RTSa}
\\
&\int_0^\infty{\rho\tilde\kappa}(1-a)(J_\nu-B_\nu(T))\dd\nu=0\,,
 \\ &
I_\nu(0,\mu) =Q_\nu(\vx,\mu), ~~ \mu<0, \quad I_\nu(H,\mu)=0,~~\mu>0.
\end{aligned}
\end{equation}
with $\Gamma$ the plane $x=0$, $Q_\nu=Q^0 B_\nu( T_s)\mu$, $\mu:=\cos\alpha$ and $\alpha$ the angle of the ray $\vom$  with $Ox$.

\section{The Method of Characteristics}

Let us derive a closed form solution of \eqref{RTa} by computing the long time solution of
\begin{equation}
\label{transport}
\partial_\tau I(\vx,\tau) +\vom\nabla_\vx I(\vx,\tau) + \kappa(\vx)I(\vx,\tau) = S(\vx).
\end{equation}
While it plays no part, notice that $\tau=t/c$.
Let $K(\vx,\vom,\tau)=\int_0^\tau\kappa(\vx+\vom s)\dd s$.  Then
\[
\frac{\dd}{\dd \tau}\left( I(\vx+\vom \tau,\tau)\e^{K(\vx,\vom,\tau)}\right)
= \e^{K(\vx,\vom,\tau)}\left[\partial_\tau I +\vom\nabla_\vx I + \kappa(\vx)I\right]_{|_{\vx+\vom \tau,\vom,\tau}}
\]
Consequently
\[
I(\vx+\vom \tau,\tau)\e^{K(\vx,\vom,\tau)} = I(\vx,0) + \int_0^\tau \e^{K(\vx,\vom,s)} S(\vx+\vom s)\dd s.
\]
Denote $\vx'=\vx+\vom \tau$; the above is also
\[
I(\vx',\tau)\e^{K(\vx'-\vom \tau,\vom,\tau)} = I(\vx'-\vom \tau,0) + \int_0^\tau \e^{K(\vx'-\vom \tau,\vom,s)} S(\vx+\vom s)\dd s.
\]
Denote by $\tau_{\vx',\vom}$ the value of $\tau$ which brings $\vx'-\vom \tau$ to the boundary:
\[
\vx'-\vom \tau_{\vx',\vom}\in\Sigma, \hbox{ and denote } \vx_\Sigma(\vx',\vom):= \vx'-\vom \tau_{\vx',\vom}.
\]
Then the long time solution of \eqref{transport} is
\begin{equation*}
\begin{aligned}
I(\vx')& = I(\vx_\Sigma(\vx',\vom))\e^{-K(\vx_\Sigma(\vx',\vom),\vom,\tau)} 
\cr&
\hskip 2cm
+ \int_0^{\tau_{\vx',\vom}} \e^{K(\vx'-\vom s,\vom,s)-K(\vx_\Sigma(\vx',\vom),\vom,\tau_{\vx',\vom})} S(\vx'-\vom (\tau_{\vx',\vom}-s))\dd s.
\end{aligned}
\end{equation*}
It can also be written as:
\begin{equation}
\begin{aligned}
I(\vx,\vom)&=I(\vx_\Sigma(\vx,\vom))\e^{-\int_0^{\tau_{\vx,\vom}}\kappa(\vx-\vom s)\dd s} + \int_0^{\tau_{\vx,\vom}} \e^{-\int_0^s\kappa(\vx-\vom s')\dd s'} S(\vx-\vom s)\dd s
\\
&=I(\vx_\Sigma(\vx,\vom))\e^{-\int_{[\vx,\vx_\Sigma]}\kappa} + \int_{[\vx,\vx_\Sigma]} \e^{-\int_{[\vx,\vx']}\kappa} S(\vx')\dd x'\,.
\end{aligned}
\end{equation}
Averaging in $\omegav$  and using \eqref{RTc}, and knowing that the solid angle integral of $f(x-\vom s)\dd\omega$ is the surface integral on $\Gamma$ of  $f(\vx'-\vx)\cdot\vn(\vx')|\vx'-\vx|^{-3}\dd\Gamma(x')$ leads to :
\begin{equation}\label{general}
\begin{aligned}
J_\nu(\vx)&:=\frac1{4\pi}\int_\SS I(\vx,\vom)\dd\omega
=\frac1{4\pi}\int_\SS I(\vx_\Sigma(\vx,\vom))\e^{-\int_0^{\tau_{\vx,\vom}}\kappa(\vx-\vom s)\dd s}\dd\omega 
\cr&
+ \frac1{4\pi}\int_\SS\int_0^{\tau_{\vx,\vom}} \e^{-\int_0^s\kappa(\vx-\vom s')\dd s'} S(\vx-\vom s)\dd s\dd\omega
\cr&
= S^E_\nu(\vx) + {\mathcal J}[S](\vx).
\end{aligned}
\end{equation}
with the notations, $z_- :=-\min(z,0)$ and 
\begin{equation}\label{sunlight}
\begin{aligned}&
 S^E_\nu(\vx):=
 =\frac{1}{4\pi}\int_\Gamma  Q_\nu(\vx,\frac{\vx'-\vx}{|\vx'-\vx|})\frac{[(\vx'-\vx)\cdot\vn(\vx')]_-}{|\vx'-\vx|^3}\e^{-\int_{[\vx,\vx']}\kappa}\dd\Gamma(\vx')
 \cr&
{\mathcal J}[S](\vx):=
 \frac1{4\pi}\int_{\Omega} S(\vx')\frac{\e^{-\int_{[\vx,\vx']}\kappa}}{|\vx'-\vx|^2}\dd x'. 
 \end{aligned}
\end{equation}
where $\vn(\vx')$ is the outer normal at $\vx'\in\p\Omega$.
The last line is true only if $\Omega$ is convex but it can be used in the general case with the following definition of  $\tilde\Omega$ and $\bar\kappa$.
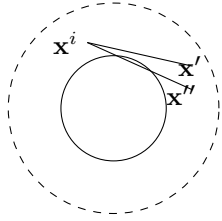
\begin{figure}[htbp]
\begin{center}
\begin{tikzpicture}[scale=0.35]
\draw (-1,2.5) node [left] {$\vx^i$};
\draw (3.5,0.5) node [left] {$\vx''$};
\draw (3.7,1.6) node [left] {$\vx'$};
\draw  (0,0) circle (2);
\draw[dashed]  (0,0) circle (4);
\draw  (3,0.7) -- (-1,2.5); 
\draw  (3,1.6) -- (-1,2.5); 
\end{tikzpicture}
\caption{Consider $\Omega$, the domain between the 2 spheres. The segment $[\vx^i,\vx']$ participates to the radiative intensity in $\vx^i$. However when the segment $[\vx^i,\vx'']$ crosses the inner sphere, the light scattered from $\vx''$ does not contribute to the light intensity in $\vx^i$.}
\label{extended}
\end{center}
\end{figure}

\begin{proposition}\label{prop:one}
Let $\bar\Omega\supset\Omega$ be a convex extension of $\Omega$. Let $\bar\kappa$ be an extension of $\kappa$ in $\bar\Omega$ such that
\[
\bar\kappa(\vx)=\kappa(vx),~~\vx\in\Omega,\quad \bar\kappa(\vx)=+\infty, ~~\vx\in\bar\Omega\backslash\Omega.
\]
Then
\[
\int_\SS\int_0^{\tau_{\vx,\vom}} \e^{-\int_0^s\kappa(\vx-\vom s')\dd s'} S(\vx-\vom s)\dd s\dd\omega
=
\int_{\bar\Omega} S(\vx')\frac{\e^{-\int_{[\vx,\vx']}\bar\kappa}}{|\vx-\vx'|^2}\dd x'.
\]
\end{proposition}\label{prop:2}
\begin{proof}: 
 If $\Omega$ is convex then $\{\vx-\omegav s ~:~s\in(0,\tau_{\vx,\omegav}), \omegav\in\SS\}=\Omega$.
Refering to figure \ref{extended}, if $\Omega$ is not convex, then the part in $\bar\Omega\backslash\Omega$ of the integral on the right is zero because $\vx'\in\bar\Omega\backslash\Omega$ implies that a portion of the integral in  the exponential will be on $[\vx^{\prime\prime},\vx']$, for some $\vx^{\prime\prime}\in\bar\Omega\backslash\Omega$. This portion  will contribute to the integral with $\bar\kappa=+\infty$ and $\e^{-\infty}=0$.
\end{proof}

\subsubsection{The Stratified Grey Case}\label{para:stratif}
The following will be helpful to assert the precision of the general algorithm \ref{algo:1} below.
  
Recall the Stefan-Blotzmann relation: $\int_0^\infty B_\nu(T)=\sigma T^4$ with $\sigma=\pi^4/15$, the Stefan constant. If $\tilde\kappa_\nu$ is constant (grey case), then \eqref{RTSa} integrated in $\nu$ yields an equation for $\bar I=\int_0^\infty I_\nu\dd\nu$:
\begin{align}
&\mu\partial_{x}\bar I +\rho\tilde\kappa\bar I=\rho\tilde\kappa \sigma T^4(x)\,,
\quad\sigma T^4(x)=\frac12\int_{-1}^1 \bar I(x,\mu)\dd\mu,
\quad x\in(0,H)\label{RTSlab3}
\end{align}
with boundary conditions $\bar I(0,\mu)|_{\mu>0}=\mu Q^0\sigma T_s^4$ and $\bar I(H,\mu)|_{\mu<0}=0$.

Even though one can find directly a solution in integral form, when $\rho\tilde\kappa$ is constant the above framework yields
\begin{equation}
\bar S^E(\vx) = \frac{Q^0\sigma T_s^4}{4\pi}\int_{\R^2} \left(\frac{\cos^2\theta}{x}\right)^2\e^{-\frac {\kappa x}{\cos\theta}}\dd\Gamma(y,z)
=\frac{Q^0\sigma T_s^4}2{\bf E}_3(\kappa x)
\end{equation}
where ${\bf E}_3$ is the third exponential integral.
Similarly  (see \cite{FGOP3})
\begin{equation}
{\cal J}[S](x) = \kappa\int_0^H  S(x'){\bf E}_1(\kappa [x'-x|)\dd x'.
\end{equation}
This approximation can be extended to the case $\kappa$ function of $x$ by a change of coordinate $x\to\tau$ with $\kappa\frac\d{\d\tau}= \frac\d{\d x}$.

\noindent
\subsection{Algorithm}\label{algo:1}
System \eqref{RTHeat} can be solved by the following iterative scheme.
\begin{enumerate}
\item {Start from $T^0\equiv 0$ and $J_\nu^0=S^E_\nu$, given by \eqref{sunlight}}
\item FOR{ $n=0,1,\dots,N-1$  }
\begin{enumerate}
\item Knowing $T^n(\vx,t)$ and ${J_\nu}^n(\vx,t)$, compute by \eqref{general}
$$
{J_\nu}^{n+1}=S_\nu^E+\mathcal J[a J_{\nu}^n+(1-a)B_\nu(T^n)],\;\nu\in(0,\infty), ~\vx\in\Omega,~t<T.
$$
\item Define $T^{n+1}$ the solution of the semilinear drift-diffusion equation
$$
\left\{
\begin{aligned}
&{}\rho c_V(\d_tT^{n+1}+\vu\cdot\nabla_x T^{n+1})-\nabla_x\cdot(\rho c_P \kappa_T\nabla_x T^{n+1})\\
&
\hskip 2cm +\mathcal B(T^{n+1})=4\pi\int_0^\infty\kappa(1-a)J_\nu^{n+1}\dd\nu\,,\\
&T^{n+1}\big\vert _{t=0}=T_{0}\,,\qquad \frac{\d T^{n+1}}{\d n}\Big\vert _{\d\Omega}=0\,,\quad\vx\in\Omega\,,\,\,t>0\,,
\end{aligned}
\right.
$$
where $\rho$ and $\vu$ are given by the Navier-Stokes equations and
\[
\ds
\mathcal B(T):=4\pi\int_0^\infty\kappa(1-a)B_\nu(\min(T_+,T_M))\dd\nu\,.
\]
\end{enumerate}
\end{enumerate}
\begin{remark}
When $\vu$ and $\kappa_T$ are negligeable, Step (b) becomes: find $T^{n+1}$ such that
\begin{equation}\label{nofluid}
\int_0^\infty\kappa(1-a)B_\nu( T_+^{n+1})\dd\nu=\int_0^\infty\kappa(1-a)J_\nu^{n+1}\dd\nu\,.
\end{equation}
\end{remark}

This equation has a unique solution because $T\mapsto B_\nu(T)$ is increasing.  As for existence, observe that the left hand-side is continuous in $T_+$, vanishes for $T_+=0$, and tends to $+\infty$ when $T_+\to+\infty$.
For a computer solution we may use a Newton method because the Hessian is positive and never vanishes.

\section{Implementation with ${\mathcal H}$-Matrices}\label{sec:5}
 Note that the computations of $S_\nu^E$ and  ${\cal J}[S]$ involve convolution like operators. Indeed,  with $\Omega$ convex bounded and  $a=0$, $S^E$ is given by \eqref{sunlight} and %
\begin{equation}
\begin{aligned}
{\mathcal J}[S_\nu](\vx) &= \tfrac1{4\pi}\int_{\Omega}S_\nu(\vx')\kappa(\vx') \frac{\e^{-\int_{[\vx,\vx']}\kappa}}{\vert \vx-\vx'\vert^2}\dd x'\,.
\end{aligned}
\end{equation}

The domain is discretized by the vertices $\{\vx^i\}_1^N$ of a tetraedral mesh. The integrals are computed with a quadrature rule using quadrature points inside the tetraedras, typically 25 points when $\vx^i$ is near to the tetraedra of the integral and 5 points otherwise ; $J$ is approximated by its $P^1$ interpolation on the mesh:
\[
J(\vx)=\sum_1^N J_j {\hat w}^j(\vx)~\text{ where ${\hat w}^j$ is the $P^1$- Finite Element hat function of vertex $\vx^j$}.
\]
 Then Step (a) of Algorithm \ref{algo:1} becomes:
\begin{equation}
\begin{aligned}&
S_{\nu,j}= a J_{\nu,j}+(1-a)B_\nu(T_j)
\\&
J_{\nu,i}=S^E_{\nu,i}+\sum_j G_{\kappa}^{ij}S_{\nu,j}\quad \text{ where } G_{\kappa}^{ij}=\frac1{4\pi}\int_\Omega \kappa\frac{\e^{-\int_{[\vx^i-\vx']}\kappa}}{|\vx^i-\vx'|^2}{\hat w}^j(\vx')\dd x'
\end{aligned}
\end{equation}

The matrix ${\bf G}$ can be compressed with the ${\mathcal H}$-matrix method so that the multiplication ${\bf G\cdot S}$ has complexity $O(N\ln N)$ for each $\nu$. The method works best when the kernel of $G$ in the integral decays exponentially with the distance between $\vx^i$ and $\vx'$. 
The ${\mathcal H}$-matrix  approximation views $\{G^{ij}\}$ as a hierarchical tree of square blocks; 
The blocks correspond to interaction between clusters of points near $\vx^j$ and near $\vx'$. Far-field interaction blocks can be approximated by a low rank matrix because their singular value decompositions (SVD) have fast decaying singular values.
 We use {\it partially pivoted adaptive cross-approximation} \cite{borm} to approach the first terms of the SVD of the blocks, because only  r-rows times r-columns columns are needed instead of the whole block, where r is the rank of the approximation.  The rank is a function of a user defined parameter $\epsilon$ in connection with the relative  Frobenius norm error. Another criteria must be met: if $R_1$ (resp $R_2$) is the radius of a cluster of points centered at $\vx_1$ (resp $\vx_2$), then one goes down the hierarchical tree  until the corresponding block  satisfies $ \max(R_1,R_2)<\eta |\vx_1-\vx_2|$ where $\eta$ is a user defined parameter. If the en of the tree is reached, the bloc is not compacted and it is displayed in red on figure \ref{hmat}.

Then, to compute the integral  $\int_0^\infty\kappa J_{\nu,i}\dd\nu$ a rectangular quadrature at points $\nu_k$ is used. 

For example when $a=0$,
\begin{equation}
\label{tournier}
\int_0^\infty\kappa J_{\nu,i}\dd\nu\approx
\sum_kS^E_{\nu_k,i}(\nu_{k}-\nu_{k-1})+\sum_k\sum_j \kappa_{\nu_k}G_{\kappa_{\nu_k}}^{ij}B_{\nu_k}({T_j})(\nu_{k}-\nu_{k-1})
\end{equation}
The same decomposition into product of $\cal H$-matrices with vectors $\{Q_{\nu_k}(\vx^i)\}_{i=1}^N$ can be applied to compute $S^E_{\nu_k}$: 
\begin{equation}
\label{tournier2}
S^E_{\nu_k,i}=\sum_1^N S^{ij}_{\kappa_{\nu_k}}Q_{\nu_k}(\vx^j), \quad
S^{ij}_{\kappa}=\int_{\Gamma} \hat w^j\e^{-\int_{[\vx^i,\vx']}\kappa}\left(\frac{[(\vx'-\vx^i)\cdot\vn(x')]^+}{|\vx'-\vx^i|^2}\right)^2\dd s(x')
\end{equation}

\subsection{Lebesgue Integrals}
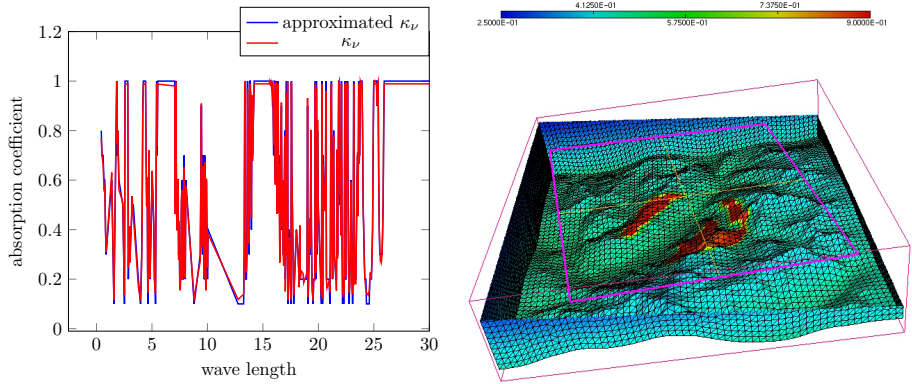
\begin{figure}[htbp]
\begin{center}
\begin{minipage} [b]{0.45\textwidth}
\begin{tikzpicture}[scale=0.7]
\begin{axis}[legend style={at={(1,1)},anchor=east}, compat=1.3,
  xmax=30, ymax=1.2,
   ylabel= {absorption coefficient},
  xlabel= {wave length}
  ]
\addplot[thick,solid,color=blue,mark=none, mark size=1pt] table [x index=0, y index=3]{fig/light.txt};
\addlegendentry{approximated $\kappa_\nu$}
\addplot[thick,solid,color=red,mark=none, mark size=1pt] table [x index=0, y index=2]{fig/light.txt};
\addlegendentry{$\kappa_\nu$}
\end{axis}
\end{tikzpicture}
\end{minipage} \hskip0.5cm
\begin{minipage}[b]{0.45\textwidth}
\includegraphics[width=6cm]{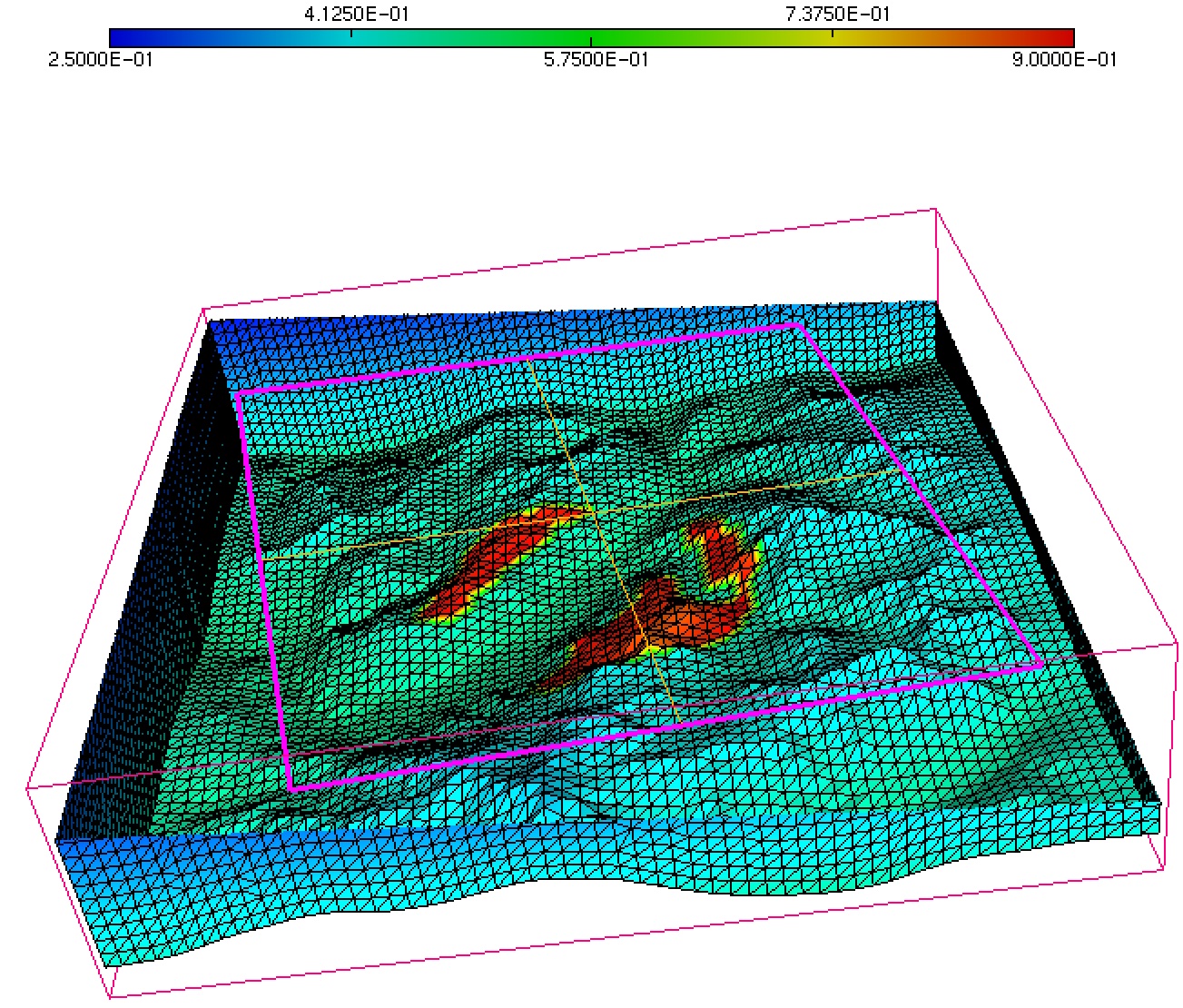}
\end{minipage}
\caption{ \label{gemini} Absorption $\kappa_\nu$ versus wavelength ($c/\nu$) in the range $(0.1,26)$ read from the Gemini data site\cite{gemini}. In the ${\mathcal H}$- step of the algorithm it is approximated by $0.01+\texttt{round}(10\kappa)/10$, shown here as ``approximated $\kappa_\nu$''.
Right: the x-dependency of $\kappa$,  the topography and the mesh. The Mont-Blanc is in the bottom left part and Chamonix is in the center. Colors are levels of the $\vx$ dependence of $\kappa$, $0.5(1-x/2)$; the red part is the intersection of the cloud with the mountains; the cloud is a cylinder circumscribing  the red parts, centered at Chamonix and between altitude 2000m and 8000m; it is used only in one simulation.}
\end{center}
\end{figure}

We will use the Gemini measurements to define $\kappa_\nu$ (figure \ref{gemini}). To represent such a function we require 683 $\nu$-points. As it is, the numerical method requires an ${\mathcal H}$-matrix for each $\nu$, but it is not feasible to define $2\times 683$ ${\mathcal H}$-matrices!

Observe that the ${\mathcal H}$-matrices ${\mathcal H,H}^E$, depend on $\nu$ directly but through $\kappa_\nu$.
It is an opportunity to reduce the number of ${\mathcal H}$-matrices down to the number of different values of $\kappa_\nu$. The the integrals with respect to $\nu$ are evaluated much like a Lebesgue integral.

\medskip

Suppose $\kappa$ takes only 2 values $\kappa_1$ and $\kappa_2$, then, according to \eqref{tournier2}, only 2 functions $S_\nu^E$ are need, $S^E_1, S^E_2$, and similarly for $G_\kappa$. Now \eqref{tournier} becomes

\begin{equation}
\begin{aligned}&
H_1^n=\sum_{\{k:\kappa_{\nu_k}=\kappa_1\}}\kappa_1(\nu_{k}-\nu_{k-1})\sum_j G_{\kappa_1}^{ij}B_{\nu_k}({T^n_j})
\\&
H_2^n= \sum_{\{k:\kappa_{\nu_k}=\kappa_2\}}\kappa_2(\nu_{k}-\nu_{k-1})\sum_j G_{\kappa_2}^{ij}B_{\nu_k}({T^n_j})
\\&
\int_0^\infty \kappa J^{n+1}_{\nu,i} = s_1\kappa_1 S^E_1 + s_2\kappa_2 S^E_2
+ H_1^n+H_2^n
\end{aligned}
\end{equation}
where $s_j$ is the measure of the $\nu$-set where $\kappa_\nu=\kappa_j$, $j=1,2.$

\medskip

There are some limitations however: $\kappa(\vx)$ must be a sum of  products of $\vx$-functions by $\nu$- functions: 
\begin{equation}\label{Nk}
\kappa_\nu(x)= \sum_j \rho_j(\vx)\kappa^j_\nu,\hbox{ where $\kappa^j_\nu$ takes only $N_k^j$ different values }\{\kappa^j_k\}_k
\end{equation}
Let $s^j_k := \{\nu~:~\kappa^j_\nu=\kappa^j_k\}.
$. Then
\begin{equation}\label{decompJ}
\begin{aligned}&
\int_0^\infty\kappa_\nu(\vx){\cal J}[S(\nu,\vx)]\dd\nu=
\cr&
\ds\frac1{4\pi}\sum_{i,j,l,k}  \rho_j(\vx)\kappa^j_k\kappa^i_l\int_\Omega\left[\int_{s^i_k}S(\nu,\vx')\dd\nu\right]\rho_i(\cdot)
\frac{\e^{-\kappa^i_l\int_{[\vx,\vx']}\rho_i(\cdot)}}{|\vx-\vx'|^2}\dd x'.
\end{aligned}
\end{equation}
Hence only the ${\mathcal H}$-matrices of kernel $\rho_i(\vx')\kappa^i_l\exp({-\kappa^i_l\int_{[\vx,\vx']}\rho_i(\cdot)})/|\vx-\vx'|^{2}$ are needed.

\smallskip

The same decomposition holds for $S_\nu^E$.
\begin{equation*}\label{decompSE}
\begin{aligned}&
\int_0^\infty \kappa_\nu(\vx) S_\nu^E(\vx)\dd\nu=
\cr&
\tfrac1{4\pi}\sum_{i,j,k,l}\rho_i(\vx)\kappa^i_l\int_{\p\Omega} \left[\int_{s^j_k}Q_\nu\dd\nu\right]\e^{-\kappa^j_k\int_{[\vx-\vx']}\rho(\cdot)}\left(\frac{[(\vx'-\vx)\cdot\vn(x')]^+}{|\vx'-\vx|^2}\right)^2\dd\Gamma(x'),
\end{aligned}
\end{equation*}
which, after discretization,  involves only the sum of products of ${\mathcal H}$-matrices of kernel on $\partial\Omega$: 
$\ds
\e^{-\kappa^j_k\int_{[\vx-\vx']}\rho_j(\cdot)}([(\vx'-\vx)\cdot\vn(x')]_-)^2/|\vx'-\vx|^4 .
$

\subsection{Computations of Exponentials on a Background Grid}

To compute $\e^{-\kappa^j_k\int_{[\vx-\vx']}\rho_j(\cdot)}$ we use a quadrature rule for the integral:
\[
\int_{[\vx-\vx']}\rho_j(\cdot) \approx \sum_{i=0}^I\rho_j(\vx_i)|\vx_{i+1}-\vx_i|,
\quad \vx_0=\vx,~~\vx_{I+1}=\vx'.
\]
To speed-up the computation of $\rho(\vx_i)$, $\rho$ is interpolated on a fine 3D Cartesian grid during the initialization phase of the computer program.

\section{Validation}

\subsection{Computing Resources}
To assert the precision and computing time we used several meshes. The grey case corresponds to $\kappa=0.5(1-x/2)$.  The
Gemini case is when $\kappa_\nu = (1-x/2)\tilde\kappa_\nu$ with $\tilde\kappa_\nu$ read from the Gemini web site.  Results are on table \ref{table:cpu} for the Chamonix valley and on figure \ref{cpu} for the flat ground case. 

\begin{table}[htp]
\caption{Computing Times in seconds (Apple M1 8 proc)}
\begin{center}
\begin{tabular}{c|c|c|c}
 nb vertices & 4053 & 30855 & 231 796\cr
\hline
Grey & 3& 19& 170\cr
Gemini & 36 & 239 & 2027\cr
\end{tabular}
\end{center}
\label{table:cpu}
\end{table}%

The code has been ported on a supercomputer. It scales perfectly for the grey case and takes 90" for 300K vertices. On the Gemini case the Newton iterations are not yet parallelized but the rest takes 360".
The memory used is governed by the compressed matrices. For the ${\mathcal H}$-matrices the compression ratio is 29 for the volumic kernel and 9 for the surfacic kernel. 

\begin{figure}[htbp]
\begin{center}
\includegraphics[width=6cm]{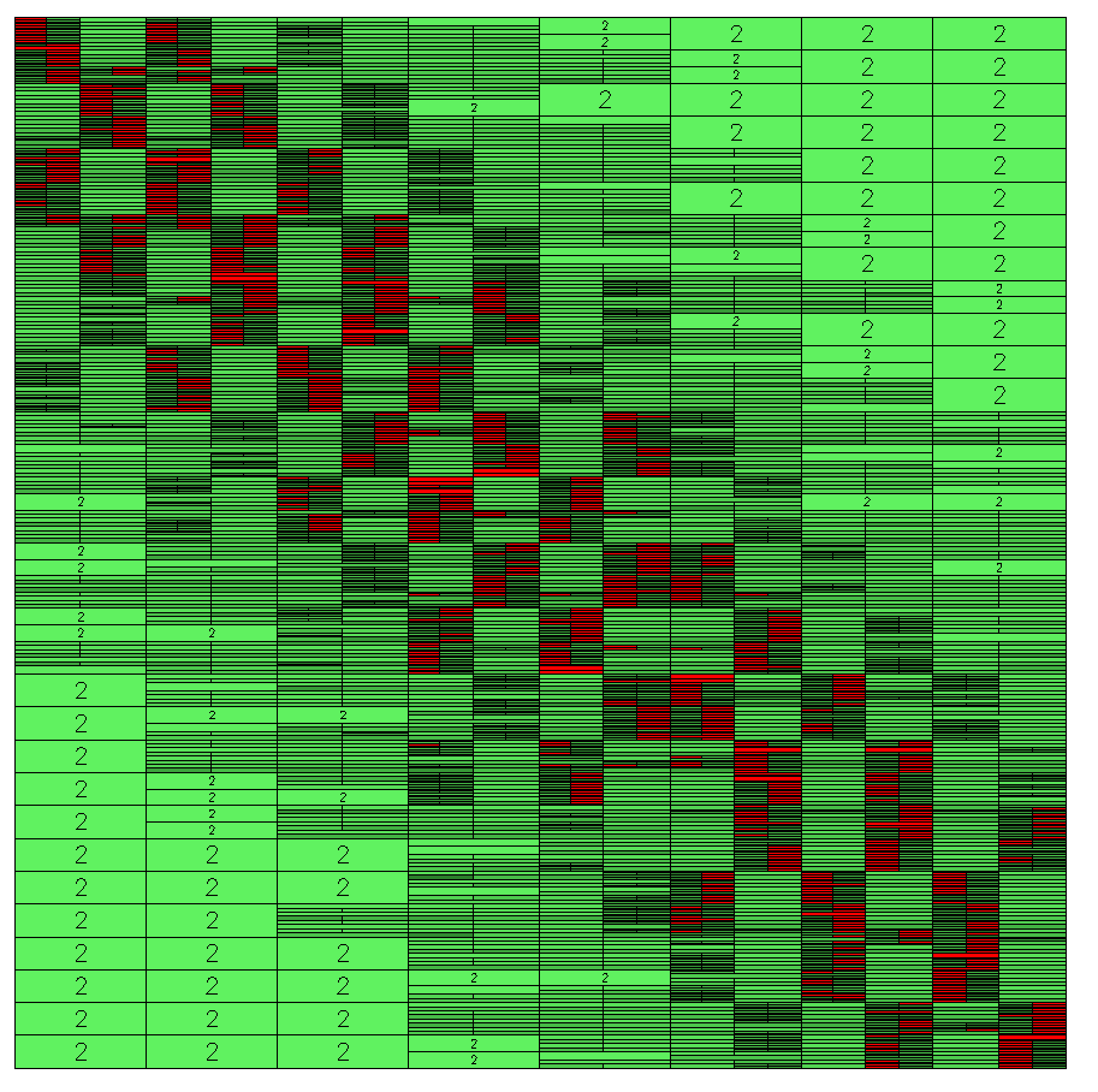}
\includegraphics[width=6cm]{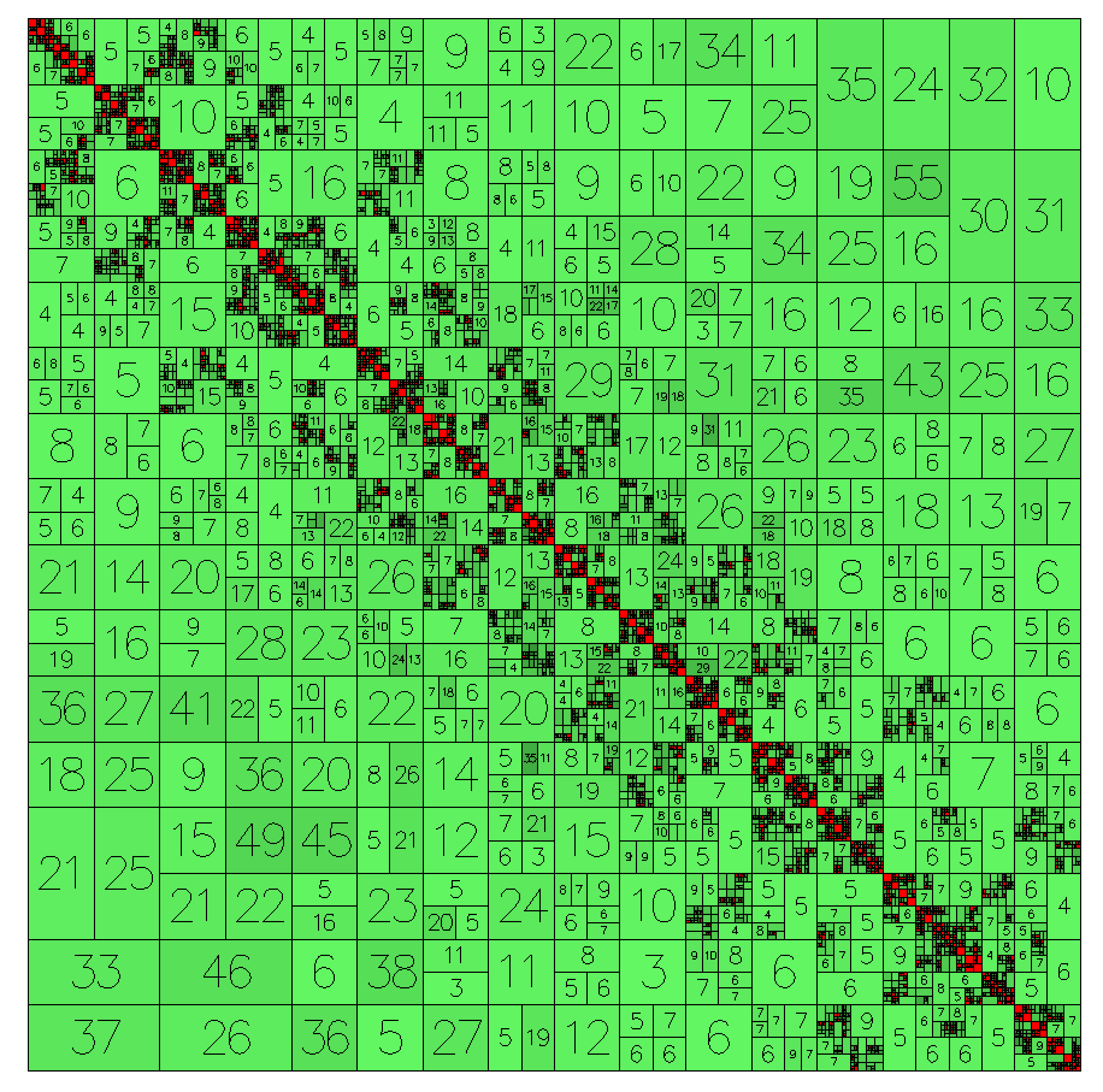}
\caption{Compressed blocks in the ${\mathcal H}$-matrices, surfacic on the left, volumic on the right. Compression is shown by colors (lighter is more compressed, red is not compressed),  numbers are the ranks of the approximations. Neighboring blocks may not have similar ranks because they maybe correspond to vertices far from each other in the physical space.}
\label{hmat}
\end{center}
\end{figure}

\subsection{Length Scale, Values for $\kappa_\nu$}
The Gemini measurements for $\kappa$ are used in a domain  $\Omega=\vx\in(0,1)\times(-0.2,3.32)\times(-3.35,0.16)$ in 10km units. The Nimbus-4 measurements \cite{nimbus} show that some infra-red radiations like $\nu_{12}=\frac3{12}\cdot 10^{14} s^{-1}$, cross the 12km thick earth atmosphere almost unaffected, which implies that $\kappa_{\nu_{12}}\sim 0$. 

Other radiations like $\nu_{15}=0.2\cdot 10^{14} s^{-1}$ are damped by 28\% (i.e. $I_{\nu_{15}}|_{x=1}=0.28 I_{\nu_{15}}|_{x=0}$, see \cite{nimbus}) which implies $\e^{-1.2\kappa_{\nu_{15}}}=0.28$, i.e. $\kappa_{\nu_{15}}=1.05$. This is associated with the presence of $\texttt{CO}_2$ in air.

 Similarly, corresponding to the presence of water in air,   $I_{\nu_{20}}$ is damped by 20\%, meaning that $\kappa_{\nu_{20}}\sim 0.3$.
 
Therefore the Gemini data can be used without scaling.  For a grey model, $\kappa=0.5$  corresponds to the average loss from the theoretical black body radiation at ground level to the actual measurements in space. Surely this value is too large in the visible  range, but the correlation between the infra-red and the visible ranges is weak.

\subsection{Temperatures on a Flat Land Exposed to Sunrays}
Our purpose here is to validate the results against the stratified numerical solutions of \cite{FGOP3}. We begin with a grey case, $\kappa=0.5$.  The domain is $\Omega=(0,1\times(-L,L)\times(-L,L)$ above the Earth surface $S:~x=0$.  The radiations of the Sun cross the atmosphere unaffected; $30\%$ is reflected and $70\%$ is absorbed and then re-emitted as a black body in all directions (Lambertian reflection) with intensity  given by\eqref{Q0S}.  Since the downward travel of the radiations is ignored,
the source of light is at $x=0$, with 
\begin{equation}\label{Q0S}
Q_\nu(\vx,\omegav)=Q^0 B_\nu(T_s)\vn\cdot\omegav, \quad Q^0=2\cdot 10^{-5},\quad T_s=1.02,~\vn=(1,0,0)^T
\end{equation}

The results depend on the mesh size $h$ and $L$.  By increasing $L$ and the number of vertices $n^3$,  convergence to the stratified case is reached (figure \ref{grey}). The convergence rate is shown on figure \ref{logerror} to be approximately $h^{1.5}$. The stratified solution is computed independently by a method sketched in Paragraph \ref{para:stratif} and detailed in  \cite{FGOP2}. 

The CPU times are plotted on figure \ref{cpu} for 5 meshes with N= 616, 4056, 30855, 80898, 231796 vertices respectively. It confirms the $N\log N$ growth.  
\smallskip

A similar exercise is done with $\kappa_\nu = 0.5(1-x/2)$. Comparison with the stratified case is possible after a change of variable in the stratified code. Results are shown on figure  \ref{greyvar}.

\begin{figure}[htbp]
\begin{minipage} [b]{0.45\textwidth}. 
\begin{center}
\begin{tikzpicture}[scale=0.7]
\begin{axis}[legend style={at={(1,1)},anchor=north east}, compat=1.3,
   xlabel= {Altitude km},
  ylabel= {Temperature C}
  ]
\addplot[thick,solid,color=pink,mark=none, mark size=1pt] table [x index=0, y index=1]{fig/milneDa.txt};
\addlegendentry{L=3, n=5}
\addplot[thick,solid,color=red,mark=none, mark size=1pt] table [x index=0, y index=1]{fig/milneDb.txt};
\addlegendentry{L=5, n=5}
\addplot[thick,dotted,color=blue,mark=none, mark size=1pt] table [x index=0, y index=1]{fig/milneDc.txt};
\addlegendentry{L=10, n=7}
\addplot[thick,solid,color=black,mark=none, mark size=1pt] table [x index=0, y index=1]{fig/milneF1D.txt};
\addlegendentry{Stratified}
\end{axis}
\end{tikzpicture}
\caption{\label{grey} Temperatures $T(x,0,0)$ versus altitude $x$ when $\kappa_\nu=\tfrac12$ for some values of $L$ and $n$. Convergence to the stratified 1D solution is visible.}
\end{center}
\end{minipage}
\hskip0.5cm
\begin{minipage} [b]{0.45\textwidth}. 
\begin{center}
\begin{tikzpicture}[scale=0.7]
\begin{axis}[legend style={at={(1,1)},anchor= north east}, compat=1.3, ymax=30,
   xlabel= {Altitude km},
  ylabel= {Temperature C}
  ]
%
\addplot[thin,solid,color=green,mark=none,mark size=1pt] table [x index=0, y index=1]{fig/milneDg.txt};
\addlegendentry{L=10, n=15}
\addplot[thin,solid,color=red,mark=none,mark size=1pt] table [x index=0, y index=1]{fig/milneDf.txt};
\addlegendentry{L=10,n=5}
\addplot[thin,dashed,color=red,mark=none,mark size=1pt] table [x index=0, y index=1]{fig/milneDe.txt};
\addlegendentry{L=10, n=10}
%
%
\addplot[thick,solid,color=black,mark=none,mark size=1pt] table [x index=0, y index=1]{fig/milneF1Dc.txt};
\addlegendentry{Stratified}
\end{axis}
\end{tikzpicture}
\end{center}
\caption{  \label{greyvar}Temperatures when $\rho\tilde\kappa=\tfrac12(1-0.5x)$: comparison with the stratified solution.}
\end{minipage}
\end{figure}

\begin{figure}[htbp]
\begin{minipage} [b]{0.45\textwidth}. 
\begin{center}
\begin{tikzpicture}[scale=0.7]
\begin{axis}[legend style={at={(1,1)},anchor=south east}, compat=1.3,
   ylabel= {CPU time (seconds)},
  xlabel= {number of vertices}
  ]
\addplot[thick,solid,color=black,mark=none, mark size=1pt] table [x index=0, y index=1]{fig/cpu.txt};
\addlegendentry{CPU}

\addplot[thick,dashed,color=pink,mark=none,mark size=1pt]expression[domain=0:2.5e5]{0.003*x};
\addlegendentry{$y=1.5 x$}

\end{axis}
\end{tikzpicture}
\caption{\label{cpu} CPU time to compute the grey case on a flat ground, indicating an $N\ln N$ growth where $N$ 
is the number of vertices in the mesh.}
\end{center}
\end{minipage}
\hskip0.5cm
\begin{minipage} [b]{0.45\textwidth}. 
\begin{center}
\begin{tikzpicture}[scale=0.7]
\begin{axis}[legend style={at={(1,1)},anchor= south east}, compat=1.3,
   xlabel= {$\log n$},
  ylabel= {$-\log\frac{\|T-T_{strat}\|}{\|T_{strat}\|}$}
  ]
\addplot[thick,solid,color=black,mark=none,mark size=1pt] table [x index=0, y index=1]{fig/logerror.txt};
\addlegendentry{Relative log error}
\addplot[thick,dashed,color=pink,mark=none,mark size=1pt]expression[domain=1.6:2.8]{1.5*x-1};
\addlegendentry{$y\sim 1.5 x$}
\end{axis}
\end{tikzpicture}
\end{center}
\caption{\label{logerror}  Log-log plot of the discretization error versus mesh size; the ``exact'' solution is the stratified one computed with many points.}
\end{minipage}
\end{figure}

\subsection{The Chamonix Valley: the Grey Case}

The domain is a portion of the atmosphere $\Omega=[h(y,z),H]\times[y_m,y_M][z_m,z_M]$ above the Earth surface $S:~x=h(y,z)$.  The radiations of the Sun cross the atmosphere unaffected; $30\%$ is reflected and $70\%$ is absorbed and re-emitted with intensity \eqref{Q0S}, except that $\vn$ is the normal of $S$. We can account for the hour of the day by rotating accordingly all normals.
\[
Q_\nu([h(y,z),y,z]^T,\omegav)=Q_0(\vx) B_\nu(T_s)\vn\cdot\omegav, \quad  T_s=1.02,
\]
We have simulated the case $Q_0=2\cdot 10^{-5}$ (no snow) and the case with snow above 2500m:
\[
Q_0(\vx) = 0.5\left(\beta + (1-\beta)\One_{x<h_{snow}}\right),~~\beta=0.3,~~
h_{snow} = 0.25,~~\vx\in S.
\]
It means that in case of snow the light source is 30\% of what it is without snow, the rest was reflected in the visible range before re-emission. Note that glaciers below 2500m are ignored.

\begin{figure}[htbp]
\begin{center}
\includegraphics[width=8.cm]{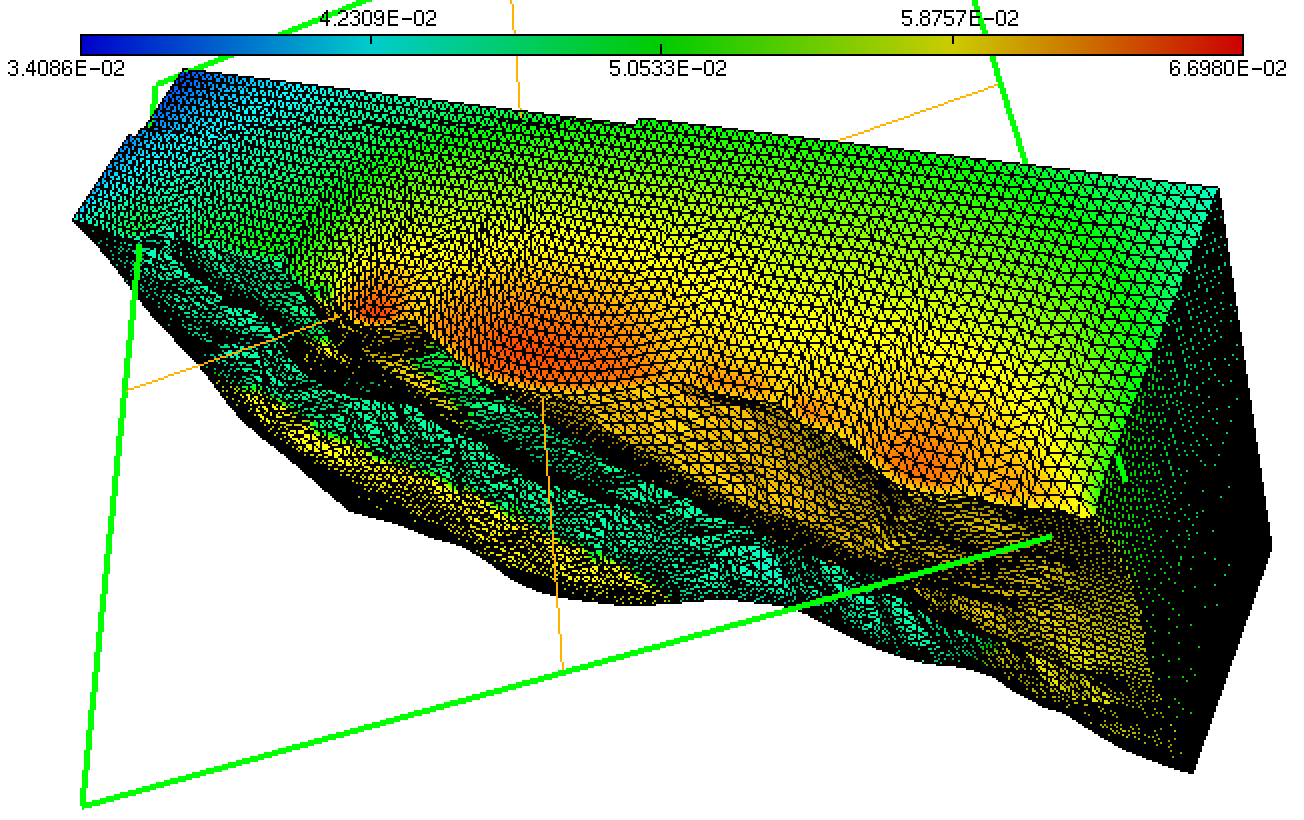}
\includegraphics[width=11.cm]{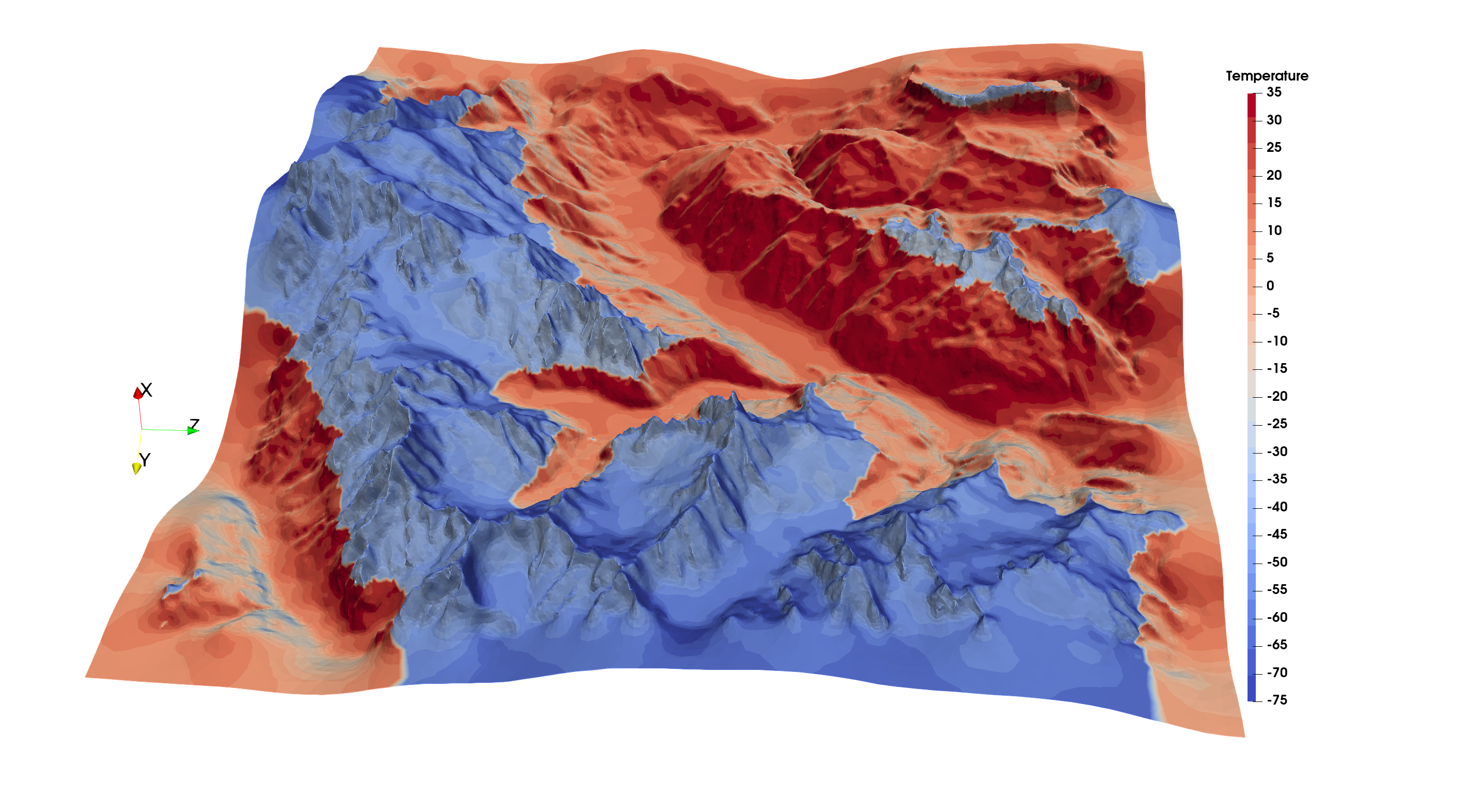}
\vskip-0.5cm
\includegraphics[width=11.cm]{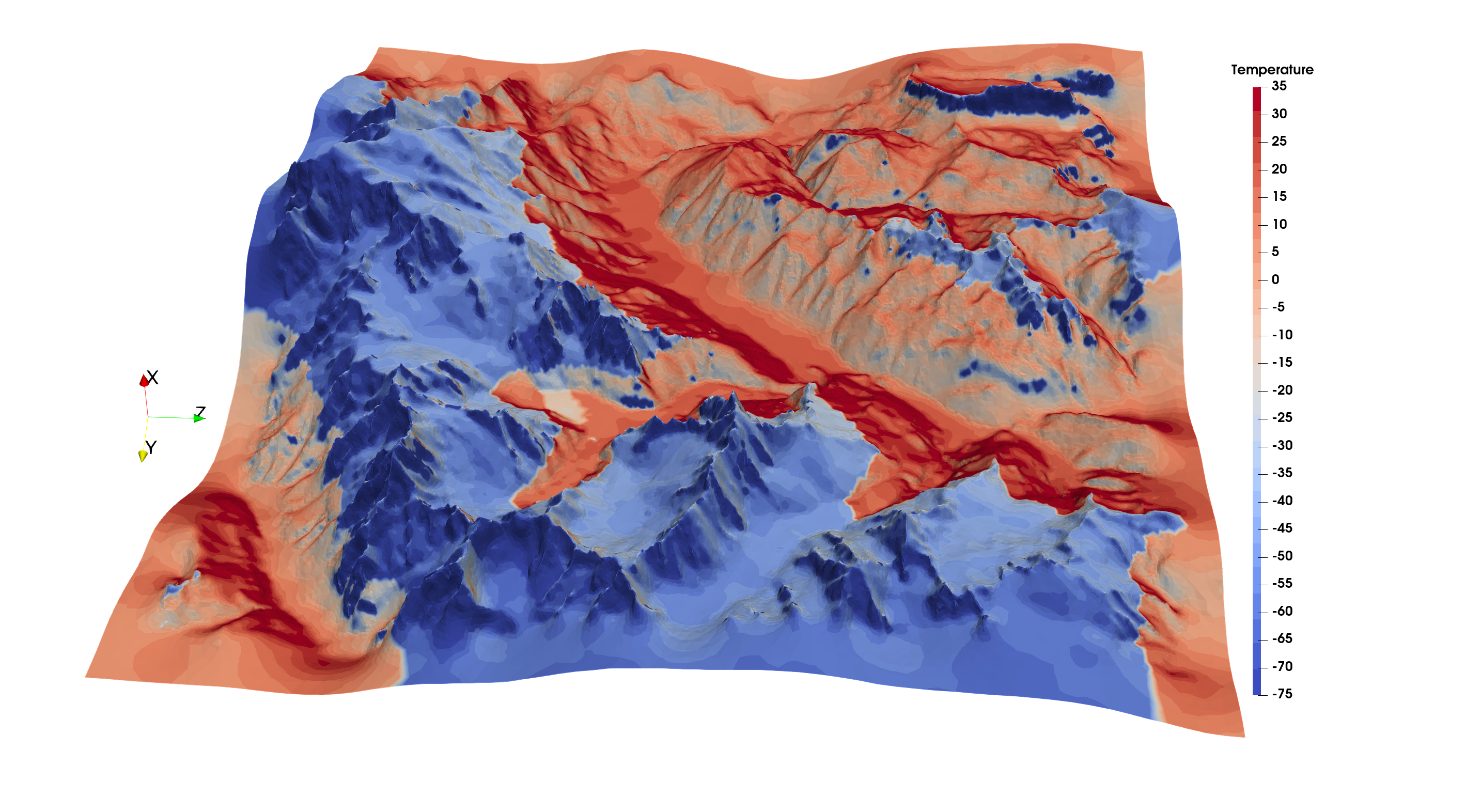}
\caption{\label{greychamhours}  All 3 computations are with $\kappa_\nu=0.5(1-x/2)$ and a clear sky (no cloud).
Top: Reduced temperatures at noon in a vertical cut of the domain (the mesh has 36K vertices).
Middle:
Ground level temperatures  in the morning: the Sun is at $45^o$ from vertical in (x,y) plane and $-20^0$ in the (x,z) plane. (Est is the lower side of the square, South is the left side). Bottom :Temperatures in the evening ((sun at $-45^o$ from vertical in (x,y) plane and $-20^0$ in the (x,z) plane; West is the upper side of the square).  Both meshes have 360K vertices.}
\label{hmat}
\end{center}
\end{figure}

As above we begin with  a $\kappa_\nu=0.5(1-x/2)$ which takes into account the rarefaction of air.  The mesh is shown on  figure \ref{greychamhours}.

Temperatures and radiations are computed in the morning and evening when the sun is at $\pm 45^o$ from the vertical and inclined towards the South by $20^0$. Results are shown on figure \ref{greychamhours}. On all such figures the altitude has been multiplied by 2 to enhance the graphics.  It can be seen that the temperature is high on the slopes exposed to the Sun.  It is seen also that the temperature is much cooler where there is snow. Temperatures are clearly too low on high mountains and perhaps too hot also on sunlit slopes. Still it is reasonable for a grey model.  
The Mont-Blanc is on the top left.  Recall that points near the border of the domain receive only half the light because of clipping, so it is unrealistically cooler near the border.

To check convergence we compared a lower solution with a higher one. To generate a lower solution initial reduced temperature in the algorithm was set to 0.01. To generate an upper solution, $T^0$ was set to 0.12.  Instead of 7 iterations, 15 were needed to reduce 
$|T^{n+1}-T^n|$ below $10^{-11}$.  Then, no visible difference could be seen between the two results.

\subsubsection{Influence of the snow}
\begin{figure}[htbp]
\begin{center}
\includegraphics[width=14cm]{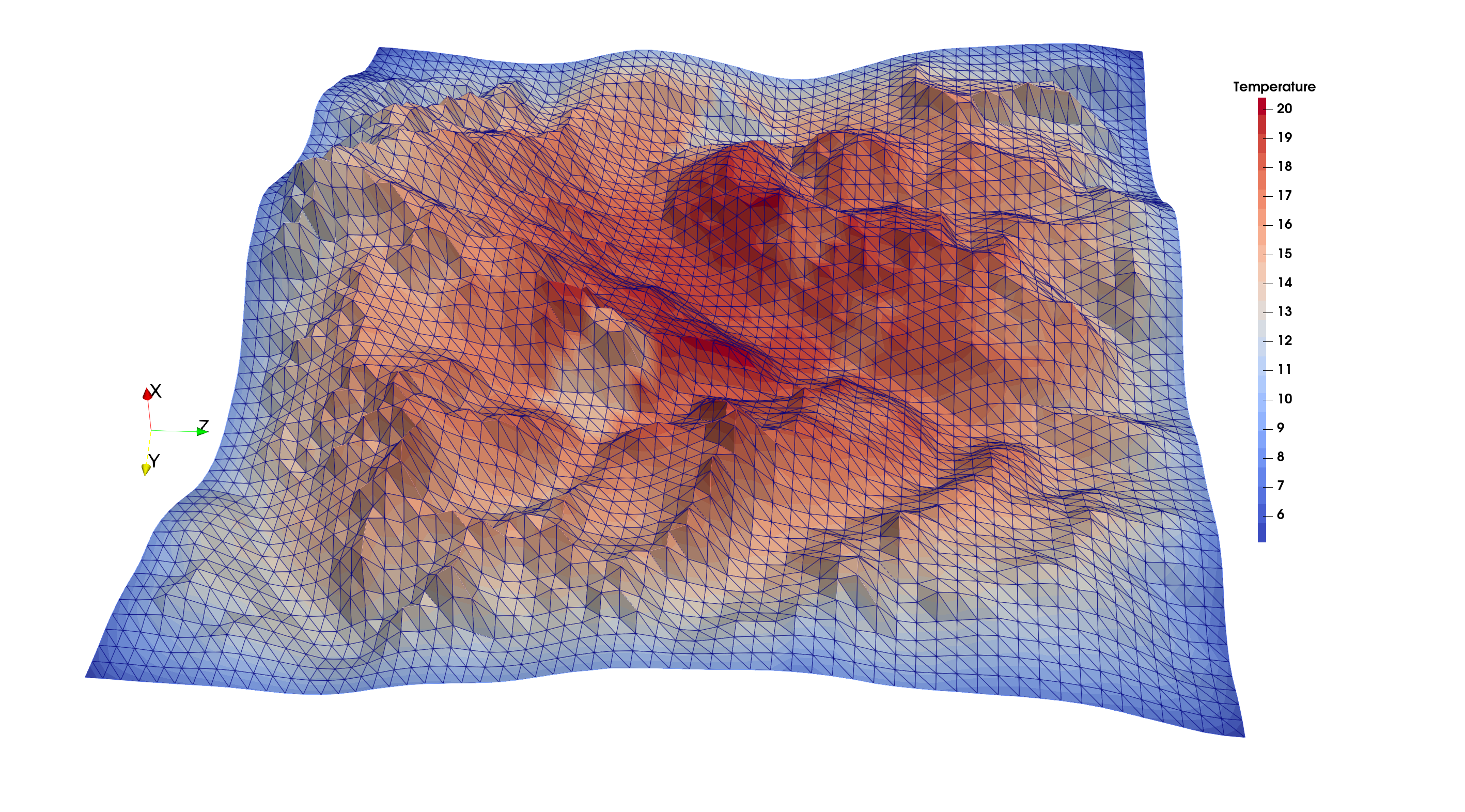}\\
\hbox{\hskip-1cm\includegraphics[width=16cm]{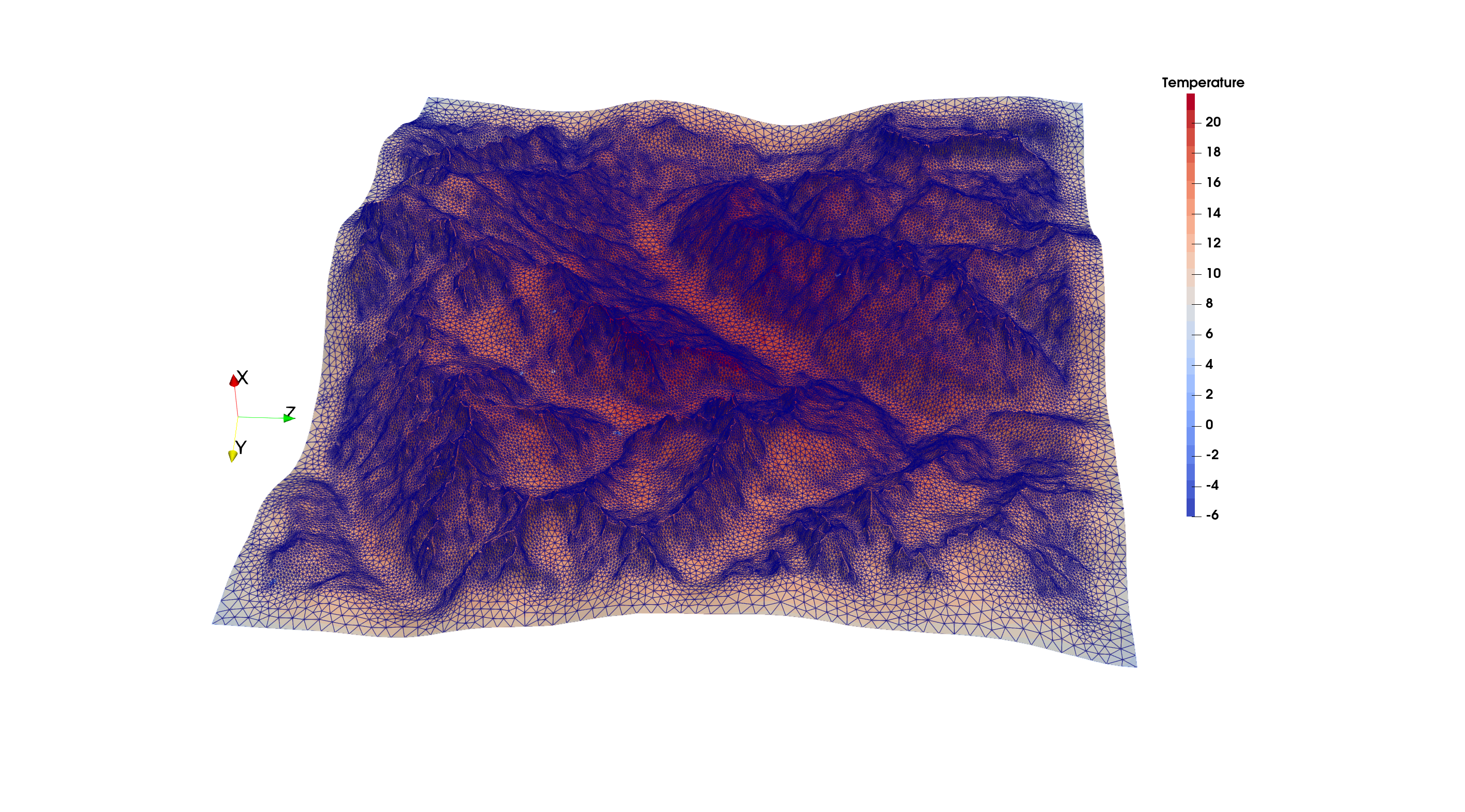}}
\caption{\label{greynosnow}Influence of the snow and the mesh: Ground temperatures  at noon  without snow. It is $2^oC$ hotter in the valley and very much hotter in the mountain, naturally. On top, the mesh has 36K vertices, and below, it has ten times more vertices.  the two pictures differ mostly near the border of the domain where the computations are wrong anyway because the points receive only one half of the sunshine other points receive (and one fourth at the corners).}
\label{hmatnosnow}
\end{center}
\end{figure}

The snow is very important for the temperature distribution at the ground level. When there is no snow 
temperatures at the ground level in the Chamonix region  at noon is $2^oC$ hotter in the valley and very much hotter in the mountain, above zero (figure \ref{hmatnosnow}).

\subsubsection{Clouds}

\begin{figure}[htbp]
\begin{center}
\includegraphics[width=12cm]{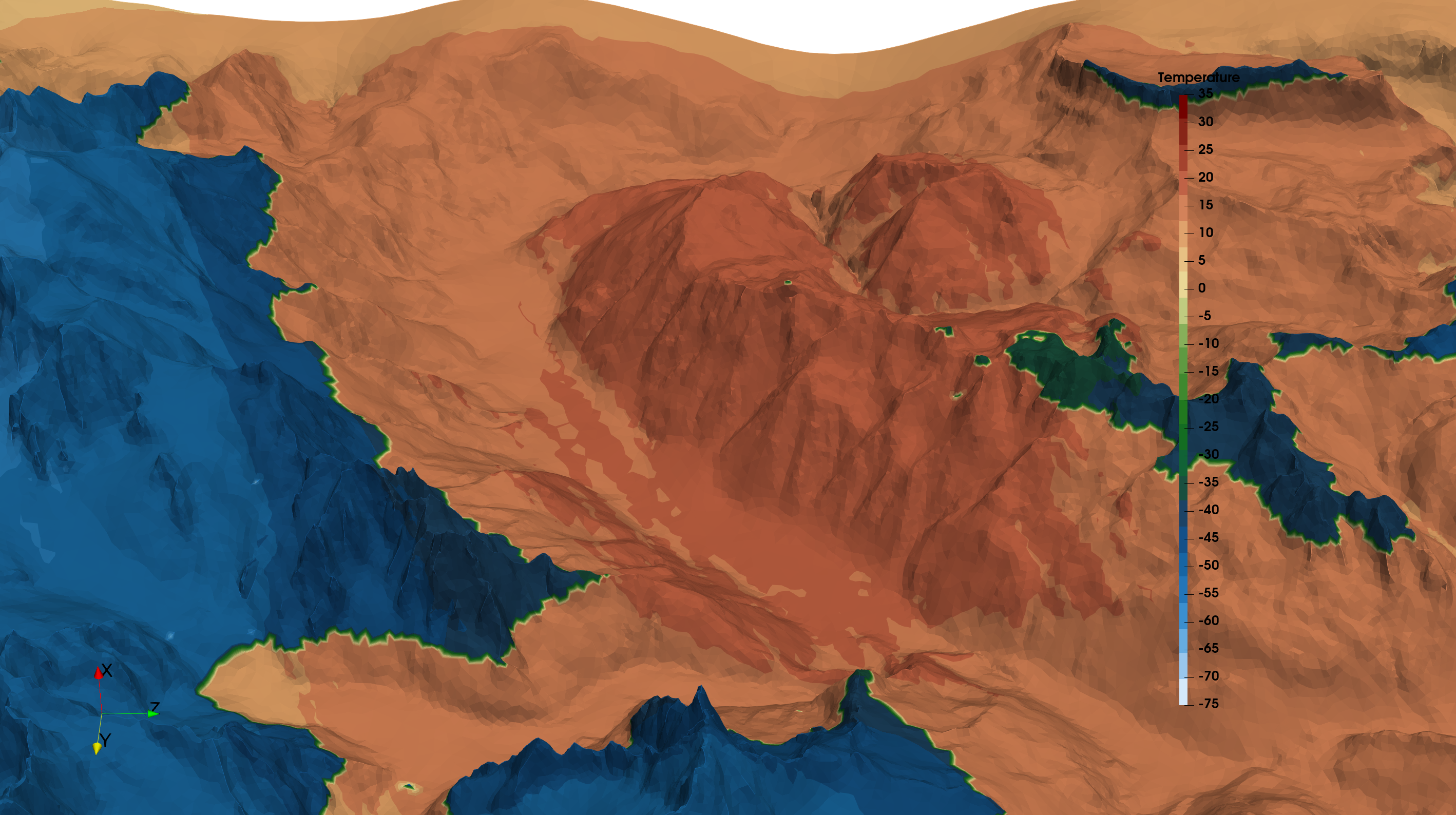}
\includegraphics[width=12cm]{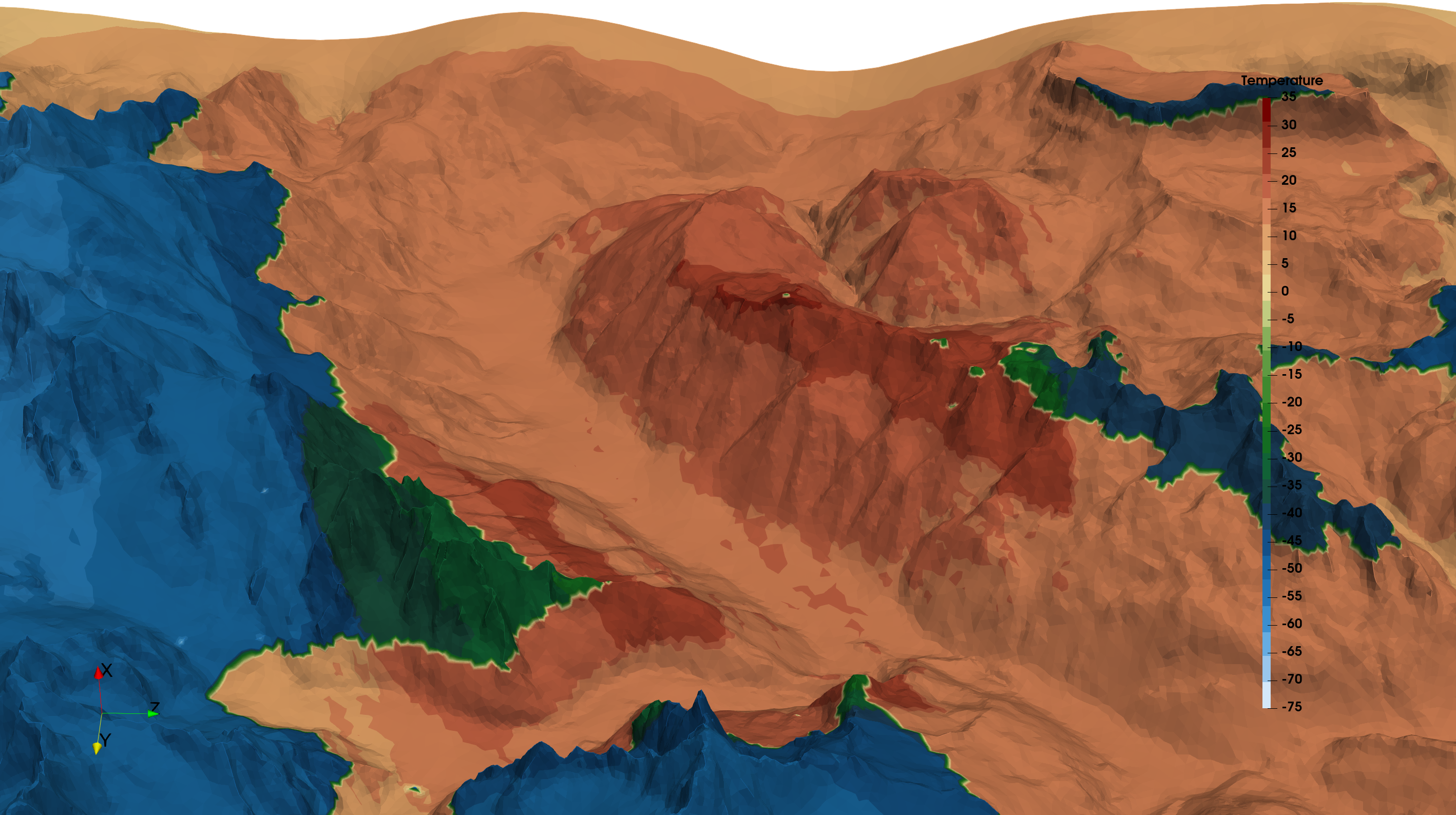}
\caption{\label{greychamnoons}Ground temperature in the Chamonix region  at noon without (top) and with (bottom) the cloud. It is hotter in the cloud but colder in the valley. }
\label{hmat}
\end{center}
\end{figure}

Let the cloud be a cylinder centered in the middle of the domain, between $x=2000$m and $x=8000$m.  There, $\kappa_\nu(\vx)$ is multiplied by $1.5$:
\[
\kappa_\nu = 0.5(1-0.5x)(1+0.5{\bf 1}[x\in(0.2,0.8):(y-1.5)^2+(z+1.5)^2<0.5]);
\]
On Figure \ref{hmat} ground temperatures are shown with and without the cloud in a zoom region with the finer mesh. It is hotter on the mountain inside the cloud but colder in the valley.

Figure \ref{tempeK} shows the temperature functions of altitude with and without the cloud (with snow) with variable density of air.

\subsection{The Chamonix Valley in the Non-Grey Case}
As in \cite{FGOP3} $\nu\mapsto \kappa_\nu$ is read from the Gemini measurements web site.

To be sure that $\int_0^\infty B_\nu(T)=\sigma T^4$ with a good precision it is necessary to extend $\tilde\kappa_\nu$ for $\nu\in(0.01,0.3$; we set it to be equal to the last available Gemini point.

Then, as shown on Figure \ref{gemini}, $\nu\mapsto \tilde\kappa_\nu$ is approximated by the nearest step function that takes only 10 values:  $\texttt{round}(10\kappa)/10$.

As explained above, by this trick we need to compute only $2\times 10$ ${\mathcal H}$-matrices, even though $\nu$-integrals are computed with 683 quadrature points.

\begin{figure}[htbp]
\begin{center}
\includegraphics[width=10cm]{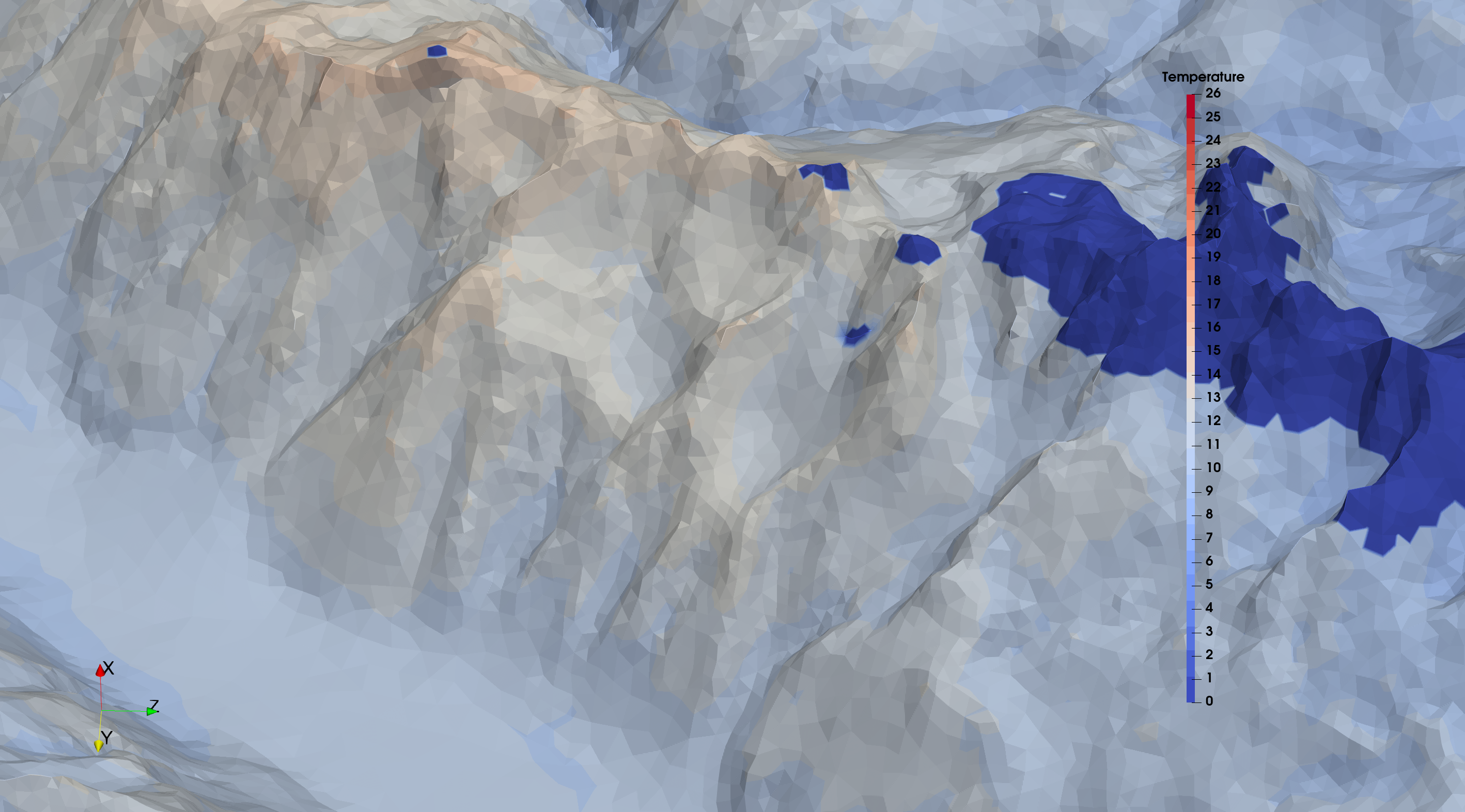}
\includegraphics[width=10cm]{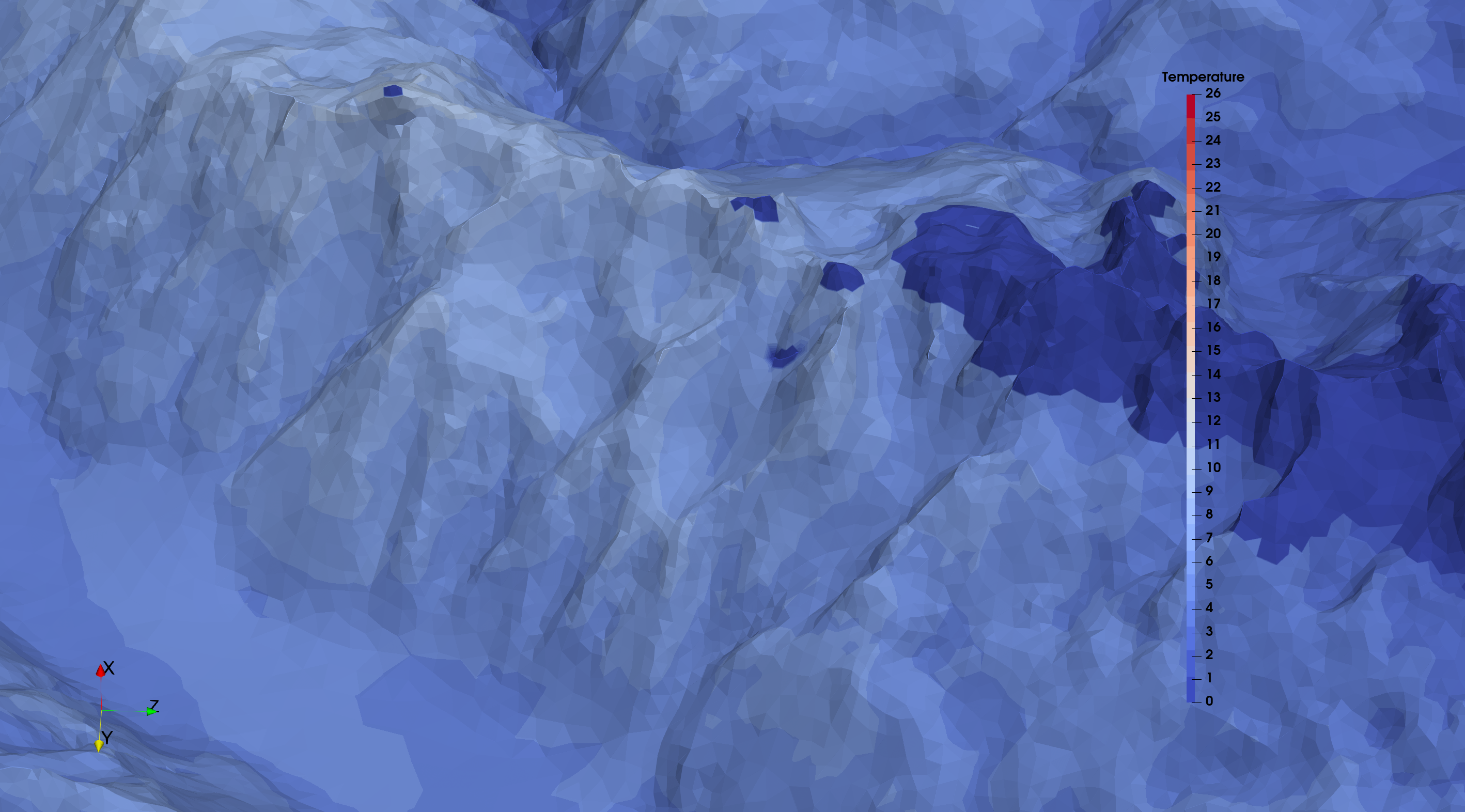}
\includegraphics[width=10cm]{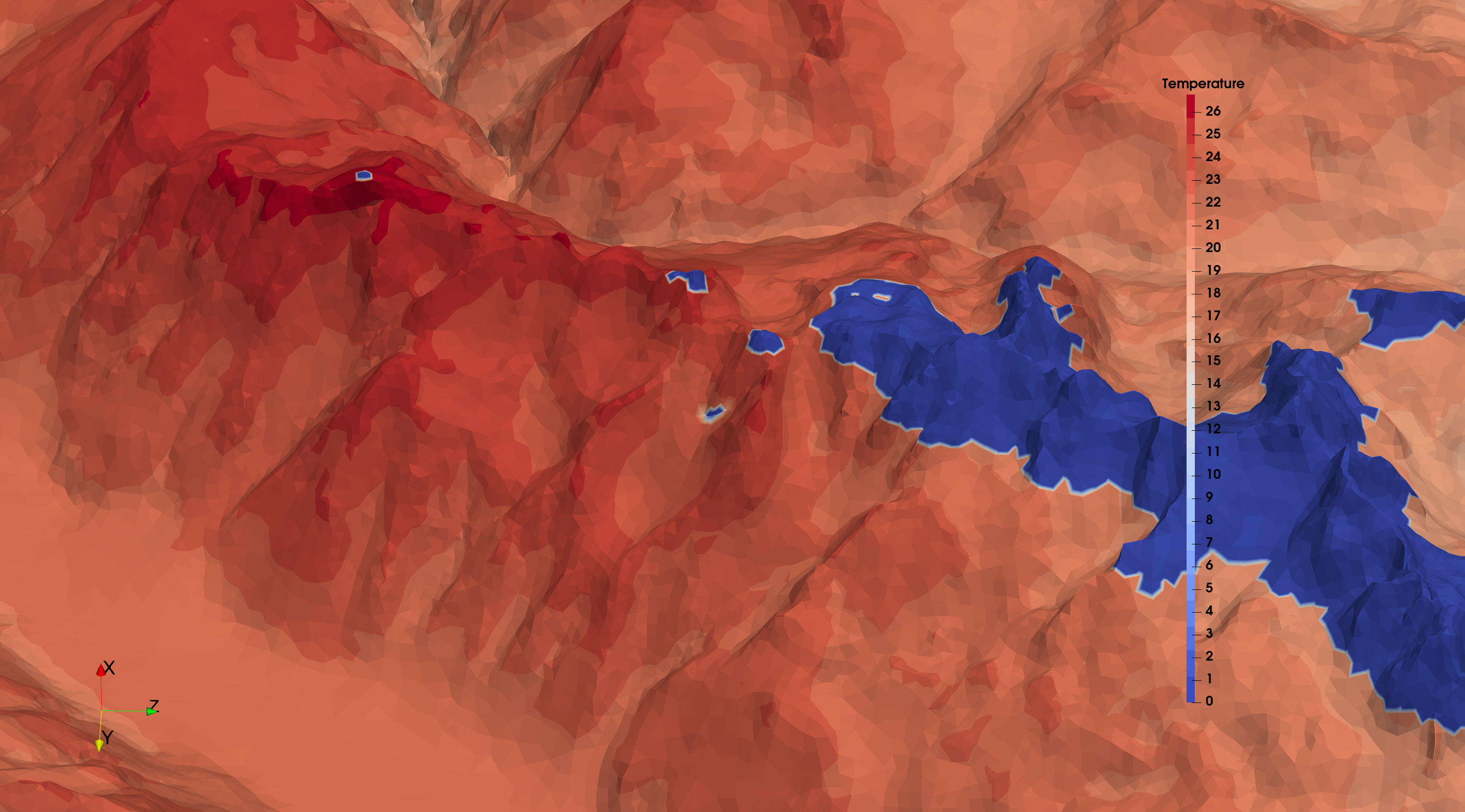}
\caption{\label{greychamnoons}Ground temperatures at noon with snow (no cloud) with $\kappa$ given by the Gemini data. Below, same situation, but the $\kappa_\nu$ from the Gemini data is put to $1-x/2$ in the range $\nu\in(c/18,c/14)$ to simulate the absorption by $\texttt{CO}_2$.  Bottom: same but with $\kappa_\nu=1-x/2$, $\nu\in(c/3,c/1.5)$ to simulate the absorption by $\texttt{CH}_4$.  The color map is the same for all 3. Bluer is colder, browner is hotter.}
\label{hmat2}
\end{center}
\end{figure}

Temperatures at noon with snow and no cloud were computed with the Gemini data of figure \ref{gemini}. Results are shown on the top picture of figure \ref{hmat2}.

Then these data were modified in the range $\nu\in(c/18,c/14)$ where $\kappa_\nu$ is set to 1.
It simulates grossly the effect of $\texttt{CO}_2$ which renders the atmosphere opaque in these frequencies.
The results are displayed on the middle picture of figure of \ref{hmat2} and on figures \ref{tempeK} and \ref{lightI}.
Note that the temperatures are colder.

A similar exercise is made but with $\kappa_\nu$  set to 1 in the range $\nu\in(c/3,c/1.5)$ to simulate an effect in the direction of  $\texttt{CH}_4$. Results are in the bottom picture of figure \ref{hmat2} and also on figures  \ref{tempeK} and \ref{lightI}. Note that the temperatures are hotter.

Hence blocking an infrared subrange can either make the day hotter or colder, depending on the position of the subrange, which in these numerical experiments we have childishly called  $\texttt{CO}_2$  and $\texttt{CH}_4$. The same observation was made in \cite{FGOP3}.

On figure \ref{lightI} the radiation intensities are shown versus wavelength ($c/\nu$) above Chamonix at altitude 3000m.

Note that the computing time is roughly ten time the one with a constant $\kappa$ because we use a Lebesgue discretization with ten levels.

\begin{figure}[htbp]
\begin{minipage} [b]{0.45\textwidth}
\begin{center}
\begin{tikzpicture}[scale=0.7]
\begin{axis}[legend style={at={(1,1)},anchor= east}, compat=1.3,
   xmin=0.1, xmax=1,
   xlabel= {Altitude km},
  ylabel= {Temperature $^o C$}
  ]
\addplot[thick,solid,color=blue,mark=none, mark size=1pt] table [x index=0, y index=1]{fig/noontempe.txt};
\addlegendentry{grey}
\addplot[thick,dashed,color=blue,mark=none, mark size=1pt] table [x index=0, y index=1]{fig/nooncloudtempe.txt};
\addlegendentry{with cloud grey}
\addplot[thick,solid,color=red,mark=none, mark size=1pt] table [x index=0, y index=1]{fig/noonKtempe.txt};
\addlegendentry{Gemini}
\addplot[thick,dotted,color=red,mark=none, mark size=1pt] table [x index=0, y index=1]{fig/noonKKtempe.txt};
\addlegendentry{Gemini+CO2}    
\addplot[thick,dotted,color=brown,mark=none, mark size=1pt] table [x index=0, y index=1]{fig/noonKKKtempe.txt};
\addlegendentry{Gemini+CH4}    
\end{axis}
\end{tikzpicture}
\caption{ \label{tempeK}Temperatures versus altitude in the center of the domain computed first with a fixed $\tilde\kappa$ with and without the cloud. Then with  a variable $\tilde\kappa$ (Gemini measurements) with $\texttt{CO}_2$ or $\texttt{CH}_4$ corrections.}
\end{center}
\end{minipage}
\hskip0.5cm
\begin{minipage} [b]{0.45\textwidth}. 
\begin{center}
\begin{tikzpicture}[scale=0.7]
\begin{axis}[legend style={at={(0,1)},anchor= west}, compat=1.3,
  xmax=30,
   ylabel= {light-intensity},
  xlabel= {wave length $\mu$m}
  ]
\addplot[thick,solid,color=red,mark=none, mark size=1pt] table [x index=0, y index=1]{fig/noonKlight.txt};
\addlegendentry{$I_\nu$ Gemini}
\addplot[thick,solid,color=blue,mark=none, mark size=1pt] table [x index=0, y index=1]{fig/noonKKlight.txt};
\addlegendentry{$I_\nu$ Gemini+CO2}
\addplot[thick,solid,color=green,mark=none, mark size=1pt] table [x index=0, y index=1]{fig/noonKKKlight.txt};
\addlegendentry{$I_\nu$ Gemini+CH4}
\end{axis}
\end{tikzpicture}
\caption{ \label{lightI} light intensity versus wave length at altitude 3000m in the center of the domain, with Gemini-$\kappa_\nu$ without and with $\texttt{CO}_2$ or $\texttt{CH}_4$ corrections.}
\end{center}
\end{minipage}
\end{figure}


\subsection{Effect of Thermal Diffusion and Wind Velocity}

\subsubsection{Wind velocity}
Convective winds due to temperature differences and modelled by a Boussinesq approximation seem out of simulation reach in areas of several kilometers.  On the other hand high atmosphere winds $u_0\sim 360km/h$ are numerically tractable if we assume that the effective viscosity is dominated by turbulence, leading to a Reynolds number around 500.  Thus the stationary Navier-Stokes equations with fixed varying density are solved in $\Omega$ for the scaled velocity $\vu/u_0$  and pressure: 
\begin{equation}\label{ns2}
\vu\cdot\n\vu + \n p -\nu\Delta\vu=0,~~\n\cdot\vu=0,\quad \vu(H,y,z)=[0,1,0]^T\quad-
\end{equation}
and Neumann conditions on the remaining boundaries $-\nu\partial_\vn\vu+p\vn=0$.

\begin{figure} [htbp]
\begin{center}
\includegraphics[width=10cm]{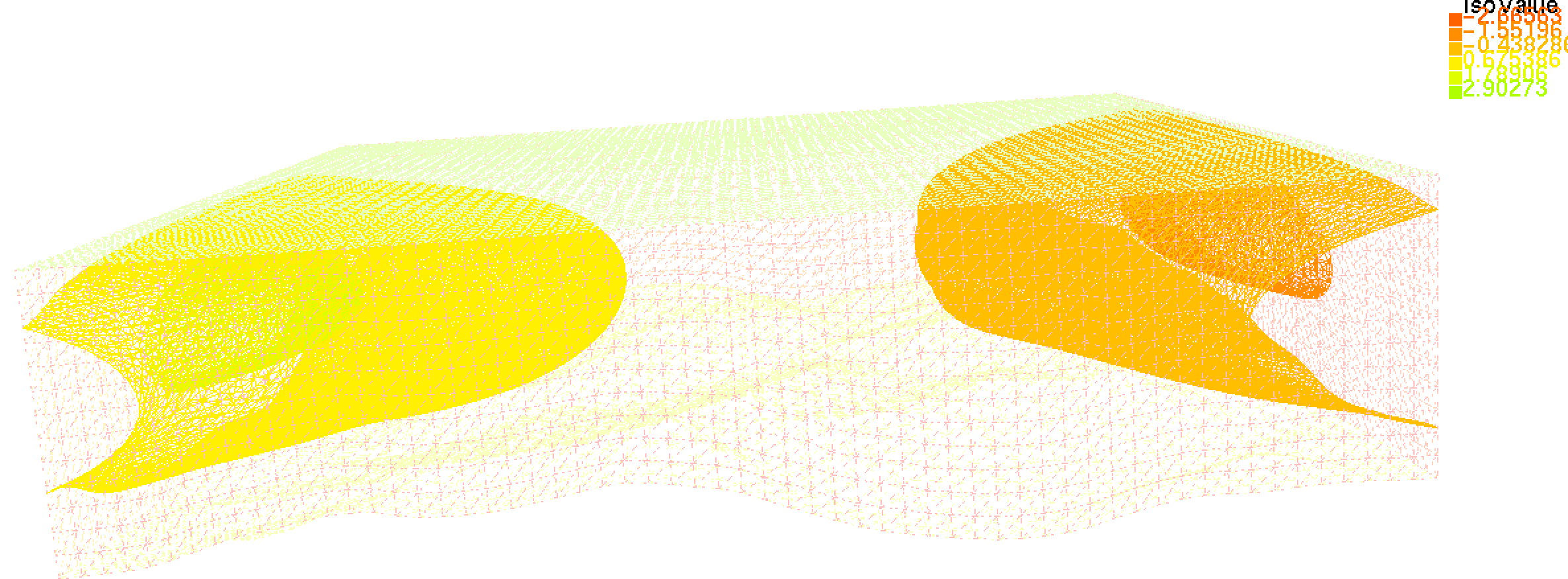}
\includegraphics[width=10cm]{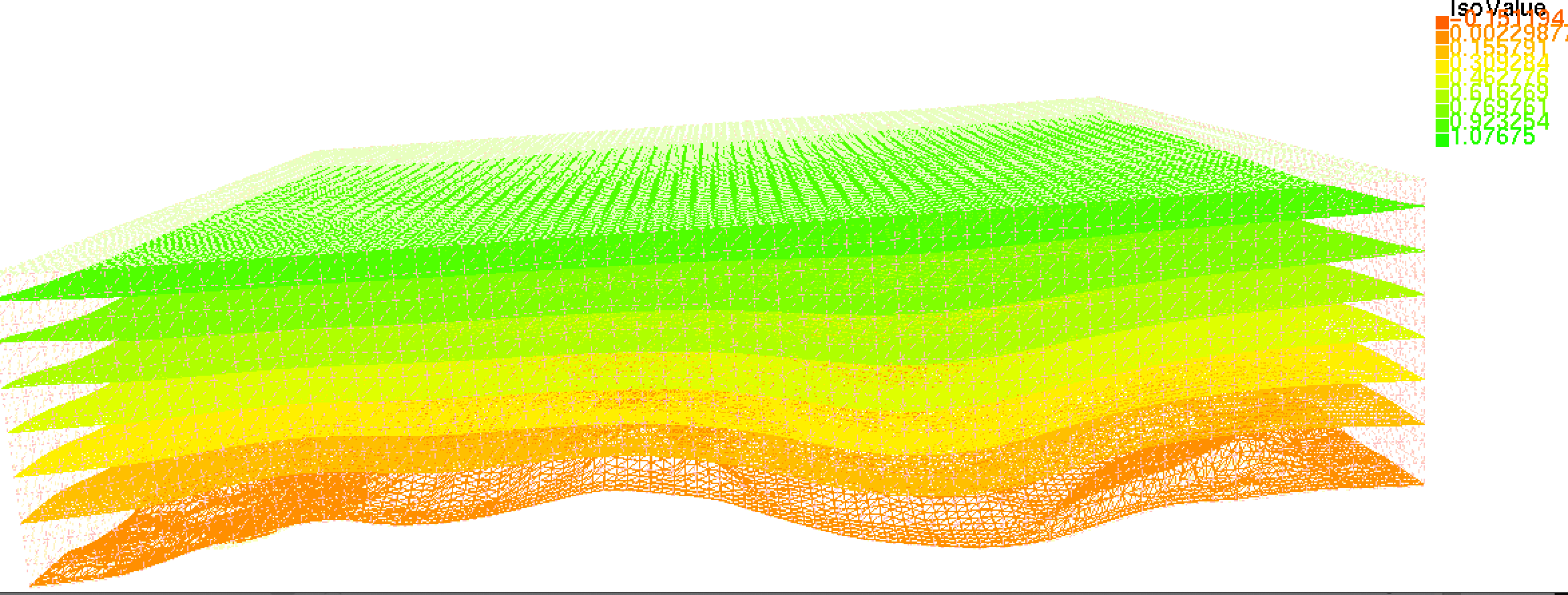}
\caption{Solution of the Navier-Stokes eqs at Re=500. Top:  $\vu_1$, Bottom: $\vu_2$
\label{fig:ns}}
\end{center}
\end{figure}

A Newton method is applied to linearize the system and the Taylor-Hood Finite Element method is used for discretization. The linear systems are solved in parallel with the MUMPS library. 
Less than a dozen iterations are sufficient for convergence in some 600 sec on the M1 (see figure \ref{fig:ns})

\
\subsubsection{Simulation of the Temperature Equation}
Recall that we have neglected the variations of $\rho$ in the diffusion term:
\[
\vu\cdot\n T -\kappa_T\Delta T = A\int_0^\infty\kappa_\nu(J_\nu-B_\nu(T))\dd\nu\hbox{ in }\Omega,\quad \frac{\partial T}{\partial n}|_{\partial\Omega}=0.
\]
In the computations, the unit length is L=10km. So $\kappa_T=2\cdot 10^{-11}$[L]$^2$/s and $\frac A{\bar\rho} = 2.70\cdot 10^{7}$[L]/s. In these units $u_0=0.01$[L]/s.  Let us divide by $u_0$:
\[
\vu\cdot\n T -\tilde\kappa_T\Delta T = \tilde A\int_0^\infty\bar\rho\kappa_\nu(J_\nu-B_\nu(T))\dd\nu\hbox{ in }\Omega,\quad \frac{\partial T}{\partial n}|_{\partial\Omega}=0,
\]
where $\vu$ is the solution of \eqref{ns2}, $\tilde\kappa_T=2\cdot 10^{-9}$, $\tilde A=5.4\cdot 10^{9}$, $\bar\rho\kappa_\nu\sim 0.5$.

\begin{figure}[htbp]
\begin{minipage} [b]{0.45\textwidth}
\begin{center}
\begin{tikzpicture}[scale=0.7]
\begin{axis}[legend style={at={(1,1)},anchor=north east}, compat=1.3,
   xmin=0.08, xmax=1,
   xlabel= {Altitude km},
  ylabel= {Temperature $^o C$}
  ]
\addplot[thick,solid,color=blue,mark=none, mark size=1pt] table [x index=0, y index=1]{fig/tempeheat.txt};
\addlegendentry{no wind}
\addplot[thick,solid,color=red,mark=none, mark size=1pt] table [x index=0, y index=1]{fig/tempeheat2.txt};
\addlegendentry{wind}

\end{axis}
\end{tikzpicture}
\caption{ \label{tempeheat}Temperatures versus altitude above Chamonix computed  with $\kappa=0.5(1-x/2)$, with snow and heat diffusion with and without wind.}
\end{center}
\end{minipage}
\hskip0.5cm
\begin{minipage} [b]{0.45\textwidth}. 
\begin{center}
\includegraphics[width=5.5cm]{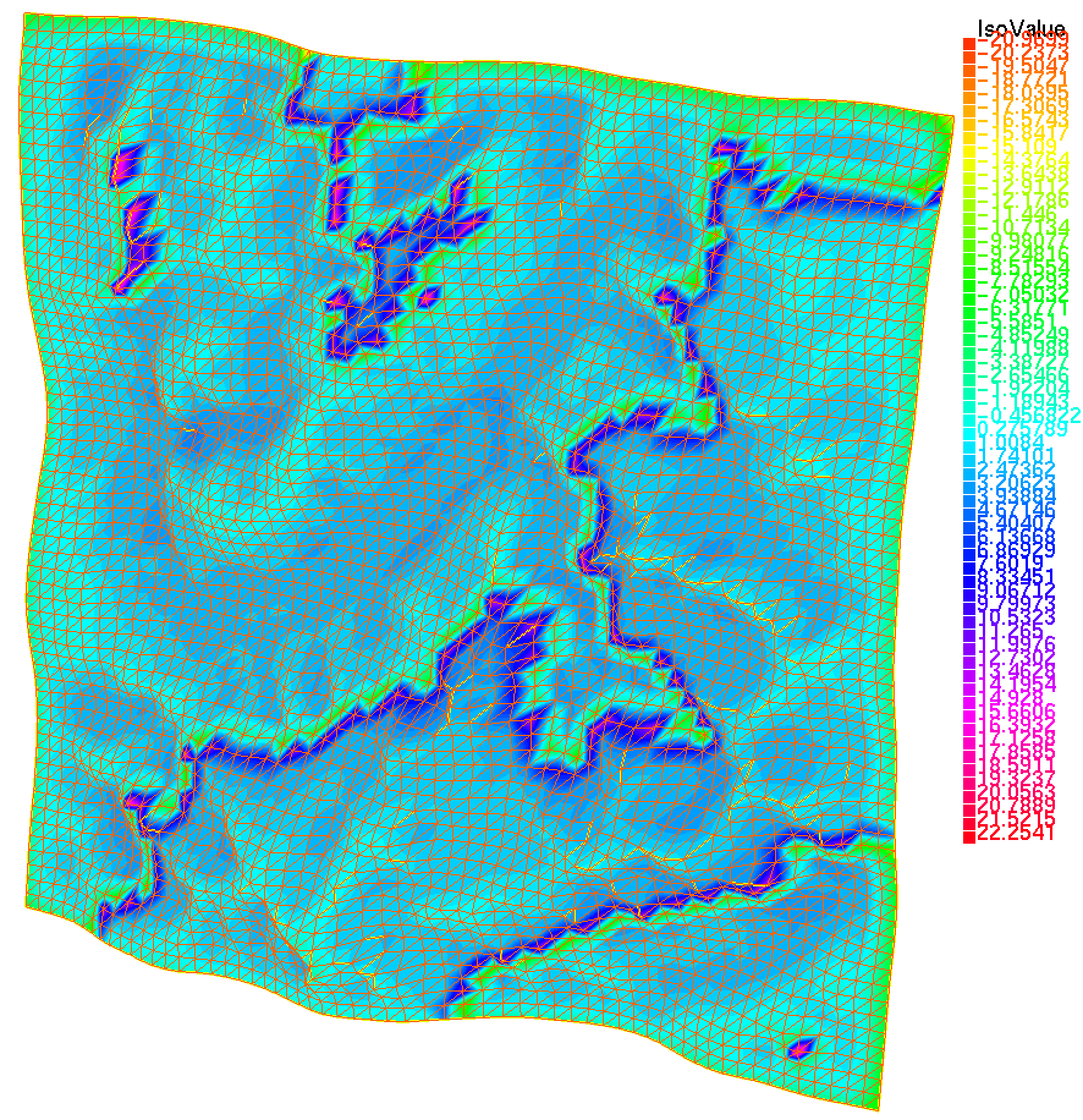}
\caption{ \label{groundheat} Temperature correction $\delta T$ on the ground due to heat diffusion, $\kappa_\nu=0.5(1-0.5x)$,  snow and a small wind.
}
\end{center}
\end{minipage}
\end{figure}

Evidently if $T^*+\delta T$ is the solution where $T^*$ is the solution with $\vu=0$ and $\tilde\kappa_T=0$, then $\delta T$ will be very small. This allows for linearization:
\[
  \vu\cdot\n \delta T -\tilde\kappa_T\Delta\delta T  + \delta T\tilde A\int_0^\infty\bar\rho\kappa_\nu B'_\nu(T^*)\dd\nu = \tilde\kappa_T\Delta\delta T^* \,.
\]
The PDE is numerically out of reach with the physical values of the parameter.  Yet to validate the concept and indicate the trend of the effect of wind and heat diffusion, we have solved it with $\tilde\kappa_T=2\cdot 10^{-3}$, $\bar\rho\kappa_\nu=(1-x/2)/2$ and $\tilde A:=A'=5.4\cdot 10^2$. Notice that if $\kappa_\nu$ does not depend on $\nu$, $A' \int_0^\infty\bar\rho\kappa_\nu  B'_\nu(T^*)=4\sigma {T^*}^3\bar\rho\kappa\sim 1.55 $.  The results are shown on figures \ref{tempeheat}, \ref{groundheat}.

\section{Conclusion}

The numerical study validates the claim that heat transfer with radiation can be solved numerically in minutes in 3D on a laptop with O(100K) mesh nodes and a fine resolution of the discontinuities of the absorption parameter in space and frequencies. On a massively parallel supercomputer with O(1M) mesh nodes it can be handled in minutes too.

The method is based on an integral formula for the radiation intensity averaged on the unit sphere of directions.  The problem is reduced to an integral equation coupled with the temperature equation with radiation as source.
An iterative scheme can be used which is mathematically shown to be convergent and monotone. The numerical errors can be estimated from the difference between the lower and the upper solutions.

A drastic numerical speed-up is obtained when the convolutions in the integrals are replaced by vector products with compressed ${\mathcal H}$-matrices. Furthermore a small number of matrices are needed only when the integrals are computed as Lebesgue integrals.

The method was tested on a portion of the atmosphere around the city of Chamonix where high mountains require a full 3D unstructured mesh.  The code is precise enough to assert the differenes due to changes in $\kappa_\nu$ is a subrange of frequenceies. All results seem physically reasonable but we make no climate claim based on the numerical results. Yet we hope to have convinced some to try the code.

\section{Acknowledgement}

Some computations were made on the machine \texttt{Joliot-Curie} of the national computing center TGCC-GENCI  under allocation A0120607330.

The computer code will soon be on github at

{\footnotesize \texttt{https://github.com/FreeFem/FreeFem-sources/tree/develop/examples/mpi/chamonix.edp}
}
%
%
%
%
%
%
%

\bibliographystyle{plain}
\bibliography{references}
\end{document}